\documentclass[11pt]{article}
\usepackage{graphicx}
\usepackage{amssymb}
\usepackage{epstopdf}
\usepackage{amsmath}
\usepackage{amsfonts}
\newcommand{\braket}[2]{\langle #1,#2 \rangle}
\newcommand{\la}{\lambda}
\newcommand{\var}{\varepsilon}

\DeclareSymbolFont{AMSb}{U}{msb}{m}{n}
\DeclareMathSymbol{\N}{\mathbin}{AMSb}{"4E}
\DeclareMathSymbol{\Z}{\mathbin}{AMSb}{"5A}
\DeclareMathSymbol{\R}{\mathbin}{AMSb}{"52}
\DeclareMathSymbol{\Q}{\mathbin}{AMSb}{"51}
\DeclareMathSymbol{\I}{\mathbin}{AMSb}{"49}
\DeclareMathSymbol{\C}{\mathbin}{AMSb}{"43}

\begin{document}
 
\addtolength{\textheight}{0 cm}
\addtolength{\hoffset}{0 cm}
\addtolength{\textwidth}{0 cm}
\addtolength{\voffset}{0 cm}

\setcounter{secnumdepth}{5}
 \newcommand{\1}{I\hspace{-1.5 mm}I}
\newtheorem{proposition}{Proposition}[section]
\newtheorem{theorem}{Theorem}[section]
\newtheorem{lemma}[theorem]{Lemma}
 \newtheorem{coro}[theorem]{Corollary}
\newtheorem{remark}[theorem]{Remark}
\newtheorem{ex}[theorem]{Example}
\newtheorem{claim}[theorem]{Claim}
\newtheorem{conj}[theorem]{Conjecture}
\newtheorem{definition}[theorem]{Definition}
 \newtheorem{application}{Application}
 
\newtheorem{corollary}[theorem]{Corollary}

\def\LX{{\cal L}(X)}
\def\LY{{\cal L}(Y)}
\def\LH{{\cal L}(H)}
 \def\ASD{{\cal L}_{\rm AD}(X)}
 \def\ASDY{{\cal L}_{\rm AD}(Y)}
\def\ASDH{{\cal L}_{\rm AD}(H)}
 \def\ASDP{{\cal L}^{+}_{\rm AD}(X)}
  \def\ASDYP{{\cal L}^{+}_{\rm AD}(Y)}
   \def\ASDHP{{\cal L}^{+}_{\rm AD}(H)}
 \def\CX{{\cal C}(X)}
\def\CY{{\cal C}(Y)}
\def\CH{{\cal C}(H)}
 \def\PX{{\cal A}(X)}
\def\PY{{\cal A}(Y)}
\def\PH{{\cal A}(H)}
\def\phi{{\varphi}}
\def\AH{A^{2}_{H}}

\title{Anti-selfdual Lagrangians II: Unbounded non self-adjoint operators and evolution equations}
\author{ Nassif  Ghoussoub\thanks{Research partially supported by a grant
from the Natural Sciences and Engineering Research Council of Canada. The
author gratefully acknowledges the hospitality and support of the Centre
de Recherches Math\'ematiques in Montr\'eal where this work was initiated.} \quad  and \quad Leo Tzou 
\\
\small Department of Mathematics,
\small University of British Columbia, \\
\small Vancouver BC Canada V6T 1Z2 \\
\small {\tt nassif@math.ubc.ca} \\
\small {\tt leo@pims.math.ca}
\\
\date{February 20, 2005}\\
}
\maketitle

\section*{Abstract} This paper is a continuation of \cite{G2}, where new variational principles were introduced  based on the concept of {\it anti-selfdual (ASD) Lagrangians}. We continue here the program of using these Lagrangians   to provide variational formulations and resolutions  to various basic equations and evolutions which do  not normally fit in the Euler-Lagrange framework. In particular, we consider stationary equations of the form  $ -Au\in \partial \varphi (u)$ as well as i dissipative evolutions of the form $-\dot{u}(t)-A_t u(t)+\omega u(t)  \in \partial \varphi (t, u(t))$ were $\phi$ is a convex potential on an infinite dimensional space.  In this paper, the emphasis is on the cases where the differential operators involved are not necessarily bounded, hence completing the results established in \cite{G2} for bounded linear operators. Our main applications deal with various nonlinear boundary value problems and parabolic initial value equations governed by the transport operator with or without a  diffusion term. 
\newpage
\tableofcontents
\section{Introduction} This paper is a continuation of \cite{G2}, where the concept of {\it anti-selfdual (ASD) Lagrangians} was shown to be inherent to many basic boundary-value and initial-value problems. A new variational framework was established where, solutions of various equations which are not normally of Euler-Lagrange type, can still be obtained as minima of functionals of the form
\[
 I(u)=L(u, A u)+\ell (b_1(x), b_2(x)) \quad  \hbox{\rm or \quad $I(u)=\int_{0}^{T}L(t, u(t), \dot u(t)+A u(t))dt+\ell (u(0), u(T)).$}
 \]
  where $L$ is an anti-self dual Lagrangian and where $A$ is  essentially a skew-adjoint operator modulo boundary terms represented by a pair of  operators $(b_1,b_2)$. For such Lagrangians,  the minimal value will always be zero and --just like the self (and antiself) dual equations of  quantum field theory (e.g. Yang-Mills and others)-- the equations associated to such minima are not derived from the fact they are critical points of the functional $I$, but because they are also zeroes of the Lagrangian $L$ itself. In other words, the solutions will satisfy
  \[
  L(u, A u)+\langle u,A u\rangle=0 \quad {\rm and} \quad L(t, u(t), \dot u(t)+A_{t}u(t))+ \langle u (t),\dot u (t) \rangle=0.
  \]
 It is also shown in \cite{G2} that ASD Lagrangians possess remarkable  permanence properties making them more  prevalent than expected and quite easy to construct and/or identify. The variational game changes from  the  analytical proofs of existence of extremals for general action functionals, to a more  algebraic search of an appropriate ASD Lagrangian for which the minimization problem is remarkably simple with value always equal to zero. This makes  them efficient new tools for proving existence and  uniqueness results for a large array of differential equations.  
 
 We tackle here again boundary value problems of the form:
 \begin{equation}
 \left\{ \begin{array}{lcl}
\label{eqn:laxmilgram1000}
 \hfill  -Au+f&\in &\partial \varphi (u)\\
\hfill  b_1(u)&=&0 \\
\end{array}\right.
\end{equation}
   as well as  parabolic evolution equations of the form:
 \begin{equation}
 \left\{ \begin{array}{lcl}
\label{eqn:1001}
 \hfill   -A u(t)-\dot{u}(t)  +\omega u (t)&\in &\partial \varphi (t, u(t)) \quad \quad \hbox {\rm a.e. $t\in [0, T]$}\\
\hfill b_{1}(u(t))&=&b_1(u_{0}) \quad \quad \quad \quad \hbox{\rm a.e $t\in [0, T]$}\\
\hfill u(0)&=&u_{0}
\end{array}\right.
\end{equation}
where $\phi$ is a convex lower semicontinuous functional and $u_{0}$ is a given initial value.
However, and unlike \cite{G2} where $A$ was assumed to be  a bounded linear operator, we deal here with existence and regularity results, hence with the more delicate framework of unbounded operators.
 
 We note that --when $\Lambda$ is linear-- such operators form a very important subset of the class of maximal monotone operators for which there is already an extensive theory (\cite{Br}, \cite{Bar}). The interest here is in the new variational approach based on the concept of anti-selfdual Lagrangians which possesses remarkable  permanence properties that maximal monotone operators either do not satisfy or do so via substantially more elaborate methods.  In a forthcoming paper (\cite{G3}), the first-named author establishes  similar results for operators of the form $
F(u)=\Lambda u +Au +\partial \varphi(u)
$
where $\Lambda$ are certain non-linear conservative operators,  $A$ are linear and positive, and $\phi$  convex, the  superposition of which is not normally covered by the theory of maximal monotone operators. 

As applications to our method, we provide a variational resolution to equations involving non self-adjoint operators such as the following {\it transport equation}:
\begin{equation}
\label{transport.eq.1}
 \left\{ \begin{array}{lcl}
    \hfill  -\vec a(x)\cdot \vec \nabla  u(x,t)+a_{0}u&=& |u(x)|^{p-1}u +f  \hbox{\rm \quad on \, $\Omega \subset \R^{n}$  
}
\\
 \hfill  u(x) &=& 0 \quad \quad \quad \quad  \quad \quad  \quad \hbox{\rm on \quad $\Sigma_{-}$, } 
       \end{array}  \right.
   \end{equation}
 where ${\vec a}=(a_{i})_{i}:\Omega \to {\bf R}^{n}$ is a vector field,  $f\in L^{2}(\Omega)$, and where $\Sigma_{-}=\{x\in \partial \Omega;\, {\vec a}(x){\hat {n}}(x)<0\}$,  ${\hat n}$ being  the outer normal vector.
 We also provide a variational resolution to general  dissipative initial value problems such as the following evolutions driven by  a superposition of the Laplacian with the transport operator.
\begin{eqnarray}
\label{Eq:transport.plus.heat}
-\frac{\partial u}{\partial t}(x,t) + \vec a(x)\cdot \vec \nabla  u(x,t) &=& \Delta_p u(x,t) + \frac{1}{2}a_0(x) u(x,t)+\omega u(x,t) \quad \hbox{\rm on $[0,T]\times \Omega$}\\
 u(x,0)& =& u_0(x)\quad \hbox{\rm on $\Omega$}\nonumber \\
 u(x,t) &=& 0\quad \quad \quad \hbox{\rm on $[0,T]\times \partial \Omega$.} \nonumber
\end{eqnarray}
But more importantly, we also deal with the more delicate case where the equation is purely non-self-adjoint such as: 
 \begin{eqnarray}
\label{Eq:transport.No.heat}
-\frac{\partial u}{\partial t}(x,t) + \vec a(x)\cdot\vec \nabla u(x,t) &=& \frac{1}{2} a_0(x)u(x,t) + u(x,t)\vert u(x,t)\vert^{p-2} + \omega u(x,t)\quad \hbox{\rm on $[0,T]\times \Omega$}\nonumber\\
u(x,0) &=& u_0(x)\quad \hbox{\rm on $\Omega$} \\
u(x,t) &= &u_0(x) \quad \hbox{\rm on $[0,T]\times \Sigma_-$}\nonumber
\end{eqnarray}
 As mentioned above, these equations  are not normally solved by the methods of the calculus of variations since they do not correspond to Euler-Lagrange equations of action functionals of the form $\int_\Omega F(x, u(x), \nabla u (x) \, dx$ or  $\int_{0}^{T}L(t, x(t), \dot x (t) dt$.
 
The paper, though sufficiently self-contained, is better read in conjunction with \cite{G2}. It is  organized as follows: In section 2, we isolate the conditions under which the composition of an anti-selfdual Lagrangian with an unbounded operator yields a Lagrangian that is also  anti-selfdual. Section 3 gives the first applications of the variational properties of ASD Lagrangians to stationary Lax-Milgram type results involving unbounded operators. In section 4 we prove the main variational principle for general Lagrangians involving semi-convex terms. 
 This principle is applied in section 5 to provide variational resolutions to several parabolic initial-value problems.
\section{ASD Lagrangians and unbounded operators}
 
We consider the class $\LX$ of  {\it convex Lagrangians} $L$ on a reflexive Banach space $X$, i.e., those functions  $L:X\times X^{*} \to \R \cup \{+\infty\}$ which are convex and lower semi-continuous (in both variables) and which are not identically $+ \infty$. The Legendre-Fenchel dual (in both variables) of $L$ is defined at any pair $(q,y)\in X^{*}\times X$ by: 
 \[
    L^*( q,y)= \sup \{ \braket{q}{x} + \braket{y}{p} - L(x,p);\, x \in X, p \in X^{*}  \}
\]
 We recall from \cite{G2} the following notions
  \begin{definition}\rm

(1) Say that  $L$ is an {\it anti-self dual Lagrangian} on $X\times X^{*}$, if  
\begin{equation}
\label{seldual}
L^*( p, x) =L(-x, -p ) \quad \hbox{\rm for all $(p,x)\in X^{*}\times X$}. 
\end{equation}
(2) $L$ is  {\it partially anti-self dual}, if 
\begin{equation}
\label{seldual.atzero}
L^*( 0, x) =L(-x, 0 ) \quad \hbox{\rm for all $x\in X$}. 
\end{equation}
 
\end{definition}
Denote by  $\ASD$ the class of anti-selfdual (ASD) Lagrangians on a given Banach space $X$. This is a quite  interesting and natural class of Lagrangians as they appear in several basic PDEs and evolution equations. The basic example of an anti-selfdual Lagrangian is given by a function $L$ on $X\times X^{*}$, of the form
\begin{equation}
L(x,p)=\varphi (x) +\varphi^{*}(-p)
\end{equation}
where $\varphi$ is a convex and lower semi-continuous function on $X$ and $\varphi^{*}$ is its Legendre conjugate on $X^{*}$.  But  the class $\ASD$ was shown in \cite{G2} to be much richer as it  goes well beyond convex functions and their conjugates, especially because it is stable under composition with skew-symmetric operators.  Indeed if $\Lambda: X\to X^{*}$ is a bounded linear skew-symmetric  (i.e., $\Lambda^{*}=-\Lambda$), and if $L$ is an ASD Lagrangian, it is then easy to see that the Lagrangian 
\begin{equation}
\label{composition}
M(x,p)=L(x, \Lambda x+p)
\end{equation}
is also anti-self dual. However, in various applications,  we are often faced with an unbounded operator $\Lambda$ which may still satisfy various aspects of anti-symmetry. In the sequel we study to what extent  the composition formula (\ref{composition}) above remains valid for such operators.

\subsection{ASD Lagrangians and unbounded skew-adjoint operators}

Let $A$ be a linear --not necessarily bounded- map from its domain $D(A)\subset X$ into $X^*$. Assuming $D(A)$ dense in $X$, we consider the domain of its adjoint $A^*$ which is defined as:
\[
D(A^*)=\{x\in X; \sup \{\langle x,Ay\rangle; y\in D(A), \|y\|_X\leq 1\}<\infty\}.
\]
\begin{definition} Let $X$ be a reflextive Banach space and let $A$ be a linear  map from its domain $D(A)\subset X$ into $X^*$. Say that
\begin{enumerate}
\item  $A$ is antisymmetric if $D(A)\subset D(A^*)$ and if $A^*=-A$ on $D(A)$.
\item  $A$ is skew-adjoint if it is antisymmetric and if $D(A)= D(A^*)$. 
\end{enumerate}
\end{definition} 
We shall also deal with situations where operators are skew-adjoint provided one takes into account certain boundary terms. We introduce the following notion

\begin{definition}  \rm 
Let $A$ be a linear map from its domain $D(A)$ in a reflexive Banach space $X$ into $X^*$ and consider $(b_1, b_2)$ to be a pair of linear maps from its domain $D(b_1, b_2)$ in $X$ into the product of two Hilbert spaces $H_1 \times H_2$. We say that $A$ is {\it skew-adjoint modulo the boundary operators $(b_1,b_2)$} if the following properties are satisfied: 
 
\begin{enumerate}
\item The set $S=D(A)\cap D(b_1, b_2)$ is dense in $X$
\item  The space $X_0:=Ker (b_1,b_2)\cap D(A)$ is dense in $X$
\item The image of $S$ by $(b_1,b_2)$ is dense in $H_1\times H_2$. 
 \item An element $y$ in $X$ belongs to $S$ if and only if 
 \[
\sup\left\{\braket{y}{Ax}- \frac{1}{2}(\| b_1(x)\|_{H_1}^2+ \| b_2(x)\|_{H_2}^2); x\in S, \| x\|_X<1\right\} <\infty.
\]
 
\item For every $x, y \in S$, we have 
 \[
\braket {y}{Ax}= -\braket{Ay}{x}+ \braket {b_1(x)}{b_1(y)}_{H_1} -\braket {b_2(x)}{b_2(y)}_{H_2}.
\]
\end{enumerate}
\end{definition}
It is clear that if $b_1$, $b_2$ are the zero operators on $X$,  then our definition coincides with the notion of skew-adjoint operator in definition 2.2.2).  
Here is our main result concerning the composition of ASD Lagrangians with non-necessarily bounded skew-adjoint operators. 
 \begin{proposition} Let $L:X\times X^*\to \R$ be an ASD Lagrangian on a reflexive Banach space $X$ such that 
 for ever $p\in X^*$, the function $x\to L(x,p)$ is bounded on the bounded sets of $X$. Let $A:D(A) \to X^*$ be  skew-adjoint modulo the boundary operators $(b_1,b_2):X\to H_1\times H_2$.  Then the Lagrangian defined by 
\begin{eqnarray*}
M(x,p)=\left\{ \begin{array}{l}
L(x,Ax+p)+\frac{{\| b_1(x)\|}_{H_1}^2}{2}+\frac{{\| b_2(x)\|}_{H_2}^2}{2}\quad 
\mbox{   if }x\in D(A)\cap D(b_1, b_2)\\
+\infty 
\, \mbox{\phantom{XXXXXXXXXXXXXXXXXXXX}if }x\notin D(A)\cap D(b_1, b_2)
\end{array}\right.
\end{eqnarray*}
is anti-self dual on $X$.
\end{proposition}

\noindent {\bf Proof:}  The idea is to use density of ${\rm Ker}(b_1,b_2)$ and the continuity of L in the first variable to split the space X in such a way that the supremum over the main term and the supremum over the boundary term are independent of each other. Indeed, 

If $\tilde x\in S:= D(A)\cap D(b_1, b_2)$, then
\begin{eqnarray*}
M^*(\tilde p,\tilde x) &=&
\sup\limits_{\stackrel{x\in S}{p\in X^*}} \left\{
   \braket{\tilde x}{p}+\braket{x}{\tilde p}
   -L(x,Ax+p)-\frac{{\| b_1(x)\|}_{H_1}^2}{2}
   -\frac{{\| b_2(x)\|}_{H_2}^2}{2}\right\} \\
\end{eqnarray*}
Substituting $y=Ax+p$, we get
\begin{eqnarray*}
M^*(\tilde p,\tilde x) &=&
\sup\limits_{\stackrel{x\in S}{y\in X^*}} \left\{
  \braket{\tilde x}{y-Ax}+\braket{x}{\tilde p} 
   -L(x,y)-\frac{{\| b_1(x)\|}_{H_1}^2}{2}
   -\frac{{\| b_2(x)\|}_{H_2}^2}{2}\right\}
\end{eqnarray*}
Since $\tilde x\in D(A)$,  we get from definition (2.3.5) that:
\begin{eqnarray*}
\braket{\tilde x}{Ax}=-\braket{x}{A\tilde x}+\braket{b_1(x)}{b_1(\tilde x)}
  -\braket{b_2(x)}{b_2(\tilde x)} .
\end{eqnarray*}
which yields
\begin{eqnarray*}
M^*(\tilde p,\tilde x) &=& \sup\limits_{\stackrel{x\in S}{y\in X^*}}
   \biggl\{ \braket{x}{A\tilde x}-\braket{b_1(x)}{b_1(\tilde x)}
   +\braket{b_2(x)}{b_2(\tilde x)} +\braket{\tilde x}{y}\\
  &&\quad \quad \quad \left. +\braket{x}{\tilde p} -L(x,y)-\frac{{\| b_1(x)\|}_{H_1}^2}{2}
   - \frac{{\| b_2(x)\|}_{H_2}^2}{2}\right\}
\end{eqnarray*}
Now for all $x_0\in {\rm Ker} (b_1,b_2)$, we obviously have 
$b_1(x)=b_1(x+x_0)$ and $b_2(x)=b_2(x+x_0), $
so that for all $x_0\in {\rm Ker} (b_1,b_2)\cap D(A)\subseteq S,$
\begin{eqnarray*}
M^*(\tilde p,\tilde x) &=&
\sup\limits_{\stackrel{x\in S}{y\in X^*}} \biggl\{\braket{x}{A\tilde x}
   -\braket{b_1(x+x_0)}{b_1(\tilde x)} +\braket{b_2(x+x_0)}{b_2(\tilde x)}
   +\braket{\tilde x}{y}\\
& &\quad \quad 
\quad +\left.\braket{x}{\tilde p}-L(x,y)-\frac{{\| b_1(x+x_0)\|}_{H_1}^2}{2}
   -\frac{{\| b_2(x+x_0)\|}_{H_2}^2}{2}\right\}
\end{eqnarray*}
It follows that 
\begin{eqnarray*}
M^*(\tilde p,\tilde x) &=&
\sup\biggl\{\braket{x}{A\tilde x+\tilde p} -\braket{b_1(x+x_0)}{b_1(\tilde x)}
   +\braket{b_2(x+x_0)}{b_2(\tilde x)} +\braket{\tilde x}{y}-L(x,y) \\
 & &\quad \quad \left. -\frac{{\| b_1(x+x_0)\|}_{H_1}^2}{2}-\frac{{\| b_2(x+x_0)\|}_{H_2}^2}{2};
    x\in S, y\in X^*,x_0\in\mbox{\rm Ker}(b_1,b_2)\cap D(A)\right\}
\end{eqnarray*}
Since $S$ is a linear space, we may set $w=x+x_0$ and write
 \begin{eqnarray*}
M^*(\tilde p,\tilde x) &=&
\sup\biggl\{\braket{w-x_0}{A\tilde x+\tilde p} 
   -\braket{b_1(w)}{b_1(\tilde x)} 
   +\braket{b_2(w)}{b_2(\tilde x)} 
   +\braket{\tilde x}{y} -L(w-x_0,y)\\
& &\quad \quad \quad -\frac{{\| b_1(w)\|}_{H_1}^2}{2} 
   -\frac{{\| b_2(w)\|}_{H_2}^2}{2};  w\in S, y\in X^*,x_0\in\mbox{Ker }(b_1,b_2)\cap D(A)\biggr\}
\end{eqnarray*}
Now, for each fixed $w\in S$ and $y\in X^*$, the supremum over $x_0\in$ ${\rm Ker} (b_1,b_2)\cap D(A)$ can be taken as a supremum over $x_0\in X$ since ${\rm Ker}(b_1,b_2)\cap D(A)$ is  dense in $X$ and all 
 terms involving $x_0$ are continuous in that variable. 
Furthermore, for each fixed $w\in S$ and $y\in X^*$, the supremum over $x_0\in X$ of the terms $w-x_0$ can be written as supremum over $v\in X$ where $v=w-x_0$. So setting $v=w-x_0$ we get
\begin{eqnarray*}
M^*(\tilde p,\tilde x) &=&
\sup\biggl\{ \braket{v}{A\tilde x+\tilde p} 
   -\braket{b_1(w)}{b_1(\tilde x)} 
   +\braket{b_2(w)}{b_2(\tilde x)} 
   +\braket{\tilde x}{y}-L(v,y)\\
& &\quad \quad \quad \quad \quad -\left.\frac{{\| b_1(w)\|}_{H_1}^2}{2} 
   -\frac{{\| b_2(w)\|}_{H_2}^2}{2}; 
    v\in X, y\in X^*, w\in S\right\}\\
 &=&   \sup_{v\in X}\sup_{y\in X^*}\left\{\braket{v}{A\tilde x +\tilde p}
   +\braket{\tilde x}{y}-L(v,y) \right\}\\
 & &\quad\quad +\sup_{w\in S}\left\{\braket{b_1(w)}{-b_1(\tilde x)} 
   +\braket{b_2(w)}{b_2(\tilde x)} -\frac{{\| b_1(w)\|}_{H_1}^2}{2} 
   -\frac{{\| b_2(w)\|}_{H_2}^2}{2}\right\}
\end{eqnarray*}
Since the range of $(b_1,b_2):S\to H_1\times H_2$ is dense 
in the $H_1\times H_2$ topology, the boundary term can be written as
 \[
  \sup_{a\in H_1}\sup_{b\in H_2}\left\{ 
   \braket{a}{-b_1(\tilde x)}+\braket{b}{b_2(\tilde x)}
   -\frac{{\| a\|}_{H_1}^2}{2} -\frac{{\| b\|}_{H_2}^2}{2}\right\}
=\frac{{\| b_1(\tilde x)\|}_{H_1}^2}{2} 
   +\frac{{\| b_2(\tilde x)\|}_{H_2}^2}{2}
   \]
while the main term is clearly equal to $ L^*(A\tilde x +\tilde p,\tilde x)=L(-\tilde x,-A\tilde x-\tilde p)$ in such a way that $M^*(p,\tilde x)=M(-\tilde x,-\tilde p)$ if $\tilde x\in D(A)\cap D(b_1, b_2)$. 
  
  Now if $\tilde x\notin S=D(A)\cap D(b_1, b_2)$ then
\begin{eqnarray*}
M^*(\tilde p,\tilde x) &=&
   \sup_{\stackrel{x\in S}{y\in X^*}}\left\{
   \braket{\tilde x}{y-Ax}+\braket{x}{\tilde p} 
   -L(x,y)-\frac{{\| b_1(x)\|}_{H_1}^2}{2}
   -\frac{{\| b_2(x)\|}_{H_2}^2}{2} \right\} \\
& &\quad \geq\sup_{\stackrel{x\in S}{{\| x\|}_X<1}}\left\{
   \braket{-\tilde x}{Ax}+\braket{x}{\tilde p} 
   -L(x,0)-\frac{{\| b_1(x)\|}_{H_2}^2}{2} 
   -\frac{{\| b_2(x)\|}_{H_2}^2}{2}.\right\} 
\end{eqnarray*}
Since by assumption $L(x,0)<C$ whenever ${\| x\|}_X<1$, we finally obtain that 
\begin{eqnarray*}
M^*(\tilde p,\tilde x)&\ge&\sup_{\stackrel{x\in S}{{\| x\|}_X<1}}\left\{
   \braket{-\tilde x}{Ax}+\braket{x}{\tilde p}-C
   -\frac{{\| b(x)\|}_{H_2}^2}{2} -\frac{{\| b_2(x)\|}_{H_2}^2}{2}\right\}\\
  &=& +\infty\\
  & = &M(-\tilde x,-\tilde p) 
  \end{eqnarray*}
since $\tilde x\notin S$ as soon as $-\tilde x\notin S$.
Therefore $ M^*(\tilde p,\tilde x)=M(-\tilde x,-\tilde p)$ for all $(\tilde x,\tilde p)\in X\times X^*$ and $M$ is an anti-selfdual Lagrangian. 
\subsection{The transport operator  $A =\vec a\cdot\vec \nabla u+\frac{(\vec \nabla \cdot \vec a)}{2} u$}

In this section we consider  the transport  operator $u\mapsto \vec a\cdot\vec \nabla u+\frac{(\vec \nabla \cdot \vec a)}{2} u$ on the space $X = L^p(\Omega)$, in conjunction with two trace operators  (restrictions) onto two appropriate subsets of $\partial \Omega$ . We show that this operator is skew-adjoint modulo the corresponding boundary operators, in the sense of Definition 2.2. These properties of the transport operator will be crucial  for the next sections where we establish existence results for stationary and evolution equations involving transport. 

Throughout this paper, we shall adopt the framework of Bardos in \cite{Ba}, and in particular all conditions that he imposes on $\Omega$ and on the smooth vector field $\vec a$ defined on a neighborhood of a $C^\infty$ bounded open set  $\Omega$  in $\R^n$.
Set  $X = L^p(\Omega)$ and define 
\[
\hbox{
 $\Sigma_{\pm} = \lbrace x \in \partial\Omega; \pm\vec a(x)\cdot\hat n(x)\ge 0\rbrace$ be the entrance and exit set of the transport operator $\vec a\cdot \vec \nabla$, }\]
 the corresponding Hilbert spaces:
\[
  H_1 = L^2(\Sigma_+ ; \vert \vec a\cdot\hat n \vert d\sigma), \quad 
  H_2 = L^2(\Sigma_- ; \vert \vec a\cdot\hat n \vert d\sigma).
  \]
as well as the boundary operators  $(b_1u, b_2u) = (u\vert_{\Sigma_+},u\vert_{\Sigma_-})$ whose domain is
\[
 D(b_1,b_2) = \lbrace u\in L^p(\Omega); (u\vert_{\Sigma_+},u\vert_{\Sigma_-})\in  H_1\times H_2\rbrace.
 \]
We shall consider the operator 
  \[
\hbox{ $Au = \vec a\cdot\vec \nabla u+\frac{(\vec \nabla \cdot \vec a)}{2} u$  with domain $D(A) = \lbrace u\in L^p(\Omega); \vec a\cdot\vec \nabla u+\frac{(\vec \nabla \cdot \vec a)}{2} u \in L^q(\Omega)\rbrace$ into $L^q(\Omega)$}.
 \]
Observe that $D(A)$ is a Banach space under the norm 
\begin{eqnarray*}
\Vert u\Vert_{D(A)} = \Vert u\Vert_p + \Vert \vec a\cdot\vec \nabla u\Vert_q
\end{eqnarray*}
and that $S:=D(A)\bigcap D(b_1)\bigcap D(b_2)$  is also a Banach 
space under the norm 
\begin{eqnarray*}
\Vert u\Vert_{S} = \Vert u\Vert_p + \Vert \vec a\cdot\vec \nabla u\Vert_q +\Vert u\vert_{\Sigma_+} \Vert_{L^2(\Sigma_+ ; \vert \vec a\cdot\hat n\vert d\sigma)}
\end{eqnarray*}
Under the assumptions listed above, $C^\infty(\bar\Omega)$ is dense in both spaces ([5]).

\begin{lemma} The operator $A$ is skew-adjoint modulo the boundary $(b_1,b_2)$ on the space $X$. 
 \end{lemma}
\noindent {\bf Proof:}
We check the five criteria of Definition 1.3. 
For  1) it suffices to note that $C^\infty(\bar\Omega) \subset S$ and $C^\infty(\bar\Omega)$ is dense in $X = L^p(\Omega)$. Similarly for 2), as we have  
 $C^\infty_0(\Omega) \subset Ker(b_1,b_2)\bigcap D(A)$, in such a way that $Ker (b_1,b_2)\bigcap D(A)$ is dense in $X$. 
Criteria (3) follows  by a simple argument with coordinate charts, as it is easy to show that for all $ (v_+, v_-) \in C^\infty_0(\Sigma_+)\times C^\infty_0(\Sigma_-)$ there exists $u \in C^\infty(\bar\Omega)$ such that $(u\vert_{\Sigma_+},u\vert_{\Sigma_-}) = (v_+, v_-)$.
The embedding of $C^{\infty}_0(\Sigma_\pm)\subset L^2(\Sigma_\pm ; \vert \vec a \cdot\hat n \vert d\sigma)$ is dense, and therefore the image of $C^\infty(\bar\Omega)$ under $(b_1,b_2)$ is dense in $H_1 \times H_2$.

For criteria 4), we need to check that, if $u\in X$, then it belongs to $S$ if and only if 
\begin{equation}
\label{belong}
 \sup\left\{\braket{u}{Av}- \frac{1}{2}(\| b_1(v)\|_{H_1}^2+\| b_1(v)\|_{H_2}^2) ; v\in S, \| v\|_X<1\right\} <\infty.
 \end{equation}
 The "if" direction follows directly from Green's theorem and the fact that $C^\infty(\bar\Omega)$ is dense in the Banach space $S$ under the norm $\Vert u\Vert_{S}$.\\
For the reverse implication, suppose that (\ref{belong}) holds, then obviously
\[
 \sup\left\{\braket{u}{Av} ; v\in C^\infty_0(\Omega),, \| v\|_X<1\right\}
 \] 
which means that  $\vec a \cdot \vec \nabla u +\frac{\vec \nabla\cdot\vec a}{2}u \in L^q(\Omega)$ in the sense of distribution, therefore $u\in D(A)$. Now to show that $u \in D(b_1)\bigcap D(b_2)$, we observe that if $u \in X$ and $u\in D(A)$, then $u \vert_{\Sigma_+} \in L^2_{loc}(\Sigma_+ ; \vert \vec a \cdot\hat n \vert d\sigma)$. To check that $u \vert_{\Sigma_+} \in L^2(\Sigma_+ ; \vert \vec a \cdot\hat n \vert d\sigma)$
a simple argument using Green's Theorem shows that (\ref{belong}) implies that 
\begin{eqnarray*}
\sup\lbrace\int_{\Sigma_+}uv \vert\vec a\cdot\hat n\vert d\sigma ; v\in C^\infty_o(\Sigma_+), \int_{\Sigma_+}\vert v\vert^2 \vert\vec a\cdot\hat n\vert d\sigma \leq 1 \rbrace <+ \infty
\end{eqnarray*}
which means that $u \vert_{\Sigma_+} \in L^2(\Sigma_+ ; \vert \vec a \cdot\hat n \vert d\sigma)$ and $u \in D(b_1)$. The same argument works for $u \vert_{\Sigma_-}$ and criterion 4) is therefore satisfied.\\
For condition 5),  note that   by Green's theorem we have  
\[
 \int_{\Omega} (({\bf a \cdot } \nabla u) v +\frac{1}{2}{\rm div \, {\bf a} })uvdx=-\int_{\Omega} (({\bf a \cdot } \nabla v) u +\frac{1}{2}{\rm div \, {\bf a} })uvdx +\int_{\partial \Omega}uv{\bf n}\cdot {\bf a}  d\sigma.
\]
for all $u,v \in C^\infty(\bar\Omega)$ and the identity on $S$ follows  since $C^\infty(\bar\Omega)$ is dense in $S$ for the norm $\Vert u\Vert_{S}$.  

\subsection{Anti-symmetric operators and ASD Lagrangians}
There are situations where anti-symmetric operators do not need to satisfy all the criteria for skew-adjointness in order to retain their composition property with ASD Lagrangians.  Here is one such setting.
\begin{lemma} Let $\phi$ be a convex proper lower  semi-continuous functional on $X$ with a symmetric domain $D(\phi )$, and let $A:D(A)\subset X\to X^*$ be an anti-symmetric operator such that $D(\phi )\subseteq D(A)$. 
Then the Lagrangian 
\begin{eqnarray*}M(x,p)=\left\{\begin{array}{ll}
\phi (x)+\phi ^*(Ax-p) &x\in D(\phi )\\
\infty &\mbox{\rm elsewhere}\end{array}\right.
\end{eqnarray*}
is anti-self dual on $X$.
\end{lemma}

\noindent {\bf Proof:} 
Let $(\tilde x,\tilde p)\in X\times X^*$ and suppose first that $\tilde x\in D(\phi )$. Then 
\begin{eqnarray*}
M^*(\tilde p,\tilde x)
&=& \sup_{x\in D(\phi )}\sup_{p\in X^*}\left\{\braket{x}{\tilde p}+\braket{\tilde x}{p}
   -\big(\phi (x)+\phi^*(Ax-p)\big)\right\}\\
&=& \sup_{x\in D(\phi )}\sup_{y\in X^*}\left\{\braket{x}{\tilde p}+\braket{\tilde x}{Ax-y}
   -\phi (x)-\phi^*(y)\right\}\\
&=& \sup_{x\in D(\phi )}\sup_{y\in X^*}\left\{\braket{x}{\tilde p}-\braket{A\tilde x}{x}
   -\braket{\tilde x}{y}-\phi (x)-\phi ^*(y)\right\}\\
&=& \phi (-\tilde x)+\phi^*(\tilde p-A\tilde x)\\
&=& M(-\tilde x,-\tilde p).
\end{eqnarray*}
If now $\tilde x\notin D(\phi )$, then 
\begin{eqnarray*}
M^*(\tilde p,\tilde x) &\ge& -\phi (0)+\sup_{p\in X^*}
   \left\{\braket{\tilde x}{p}-\phi^*(-p) \right\}\\
&=& -\phi (0)+\phi (-\tilde x)=+\infty= M(-\tilde x,-\tilde p).
\end{eqnarray*}

\section{Variational resolution of equations of the form $-Au\in \partial \phi (u)$}
If $L$ is in $\ASD$, then necessarily
 \begin{equation}
L(x, p)\geq -\langle x, p\rangle \quad \hbox{\rm  for every $(x, p) \in X\times X^{*}$}, 
\end{equation}
which means that 
\begin{equation}
I(x)=L(x,0)\geq 0 \quad \hbox{\rm  for every $x \in X$}, 
\end{equation}
 ASD Lagrangian are variationally interesting because the minima of $L(x,0)$ are often equal to zero as the following proposition --established in \cite{G2}-- indicates.
 
\begin{proposition} Let $L$ be a convex lower-semi continuous functional on a reflexive Banach space $X\times X^{*}$. Assume that $L$ is  a partially anti-self dual Lagrangian and that for some $x_{0}\in X$, the function  $p\to L(x_{0},p)$ is bounded above on  a neighborhood of the origin in $X^{*}$.
 Then there exists $\bar x\in X$, such that: 
 \begin{equation}
 \left\{ \begin{array}{lcl}
\label{eqn:existence}
  L( \bar x, 0)&=&\inf\limits_{x\in X} L(x,0)=0.\\
 \hfill (0, -\bar x) &\in & \partial L (\bar x,0).
\end{array}\right.
 \end{equation}
\end{proposition}
This result was used in \cite{G2} to establish variationally various existence results for operator equations which are not normally of Euler-Lagrange type. We can now deal with cases where these operator are not necessarily bounded. 

\subsection{A Lax-Milgram type result for unbounded operators}

\begin{proposition}
Let $\phi :X\to \R\cup \{+\infty\}$ be proper convex and lower semi-continuous and assume that  
$A: D(A)\subset X \to X^*$ is a skew-adjoint operator modulo the boundary  $(b_1, b_2):D(b_1,b_2)\to H_1\times H_2$ where $H_1$, $H_2$ are two Hilbert spaces. 
 Suppose there exists a constant $C>0$ such that for every $x\in X$, 
\begin{equation}
\frac{1}{C}\left(\| x\|_X^{p_1}-1\right)\leq\phi (x) \leq C\left( \| x\|_X^{p_2}+1\right),  
\end{equation}
where $p_1,p_2>0$. Then  there exists $\overline x\in D(A)\cap D(b_1, b_2)$ such that
\begin{equation}
 \phi (\overline x)+\phi^* (-A\overline x)+\frac{{\| b_1(\overline x)\|}_{H_1}^2}{2}
+\frac{{\| b_2(\overline x)\|}_{H_2}^2}{2}=0,
\end{equation}
\begin{equation}
 -A\overline x\in\partial\phi (\overline x) \quad {\rm and} \quad b_2(\overline x)=0.  
\end{equation}
 
\end{proposition}

\noindent {\bf Proof:} 
The Lagrangian $M(x,p):X\times X^*\to\overline{\R}$ defined by
\[
M(x,p) =\left\{\begin{array}{l}
\phi (x)+\phi ^*(-Ax-p)+\frac{{\| b_1(x)\|}_{H_1}^2}{2}  +\frac{{\| b_2(x)\|}_{H_2}^2}{2} \quad {\rm if}\, x\in D(A)\cap D(b_1, b_2)\\
+\infty  \hfill \text{otherwise}
\end{array}\right.
\]
is anti-selfdual on $X$. Indeed, $L(x,p):=\phi (x)+\phi ^* (-p)$  is clearly ASD, and since $x\mapsto L(x,p)$ is bounded on the bounded sets of $X$ for all $p\in X^*$,  we apply Proposition 2.1 to conclude that $M$ is ASD. Moreover, since  $p\to M(0,p)$ is bounded on the bounded sets of $X^*$, Proposition 3.1 applies and we obtain  
$\overline x$ such that $0=\inf_{x\in X} M(x,0)=M(\overline x,0)$. Since $M(\overline x,0)<\infty$, we get that   $\overline x\in D(A)\cap D(b_1, b_2)$ and $
 \phi (\overline x)+\phi^* (-A\overline{x})
  +\frac{{\| b_1(\overline x)\|}_{H_1}^2}{2} +\frac{{\| b_2(\overline x)\|}_{H_2}^2}{2}=0.$

Now observe that
\begin{eqnarray*}
0&=&\phi (\overline x)+\phi^* (-A\overline x)
  +\frac{{\| b_1(\overline x)\|}_{H_1}^2}{2} 
  +\frac{{\| b_2(\overline x)\|}_{H_2}^2}{2}\\
&\geq&-\braket{\overline x}{A\overline x}
  +\frac{{\| b_1(\overline x)\|}_{H_1}^2}{2}
  +\frac{{\| b_2(\overline x)\|}_{H_2}^2}{2}\\
&=& {\| b_2(\overline x)\|}_{H_2}^2 \ge 0, 
\end{eqnarray*}
and therefore $-A\overline x\in \partial\phi (\overline x)$ and $b_2(\overline x)=0$.

\subsection{Applications to PDE involving the transport operator}
 
We will deal with two types of equations:
\begin{enumerate}
\item Transport equation:
\begin{eqnarray*}
\vec a(x)\cdot\vec \nabla u(x) + a_0(x) u(x) + u(x)\vert u(x)\vert^{p-2} &=& f(x)\quad \mbox{for $x\in\Omega$}\\
u(x) &=&0 \quad \quad \quad \mbox{for  $x\in\Sigma_+$}
\end{eqnarray*}
\item  Transport equation with viscosity:
\begin{eqnarray*}
-\Delta_p u(x) +\vec a(x)\cdot\vec \nabla u(x) + a_0(x) u(x) + u\vert u(x)\vert^{m-2} &=& f(x)\quad \mbox{ $x\in\Omega$}\\
u(x)& =& 0 \quad \mbox{ $x\in\partial \Omega$}
\end{eqnarray*}
\end{enumerate}
We shall see that the first order differential operators $\vec a(x)\cdot\vec \nabla$ in the transport equation is skew adjoint modulo the boundary,  while in the case involving the p-Laplacian,  it will just be an anti-symmetric operator.  \\

\noindent{\bf 1. Transport Equations}\\
In this case we will assume that the domain $\Omega$ and the vector field $\vec a(\cdot)$ satisfies all the assumption in section 2.2.  We distinguish two cases:\\

\noindent   {\bf Case 1: $p\geq 2$.}\\
The Banach space is then $X = L^p(\Omega)$ since by lemma 2.4 the operator $A : D(A)\to X^*$ defined by 
$A : u\mapsto\vec{a}\cdot\vec \nabla u+\frac{1}{2}\left( \vec \nabla\cdot\vec{a}\right) u$ with domain
$D(A) = \lbrace u\in L^p(\Omega) \mid \vec a\cdot\vec \nabla u+\frac{\vec \nabla\vec a}{2} u\in L^q(\Omega)\rbrace
$
is then skew-adjoint modulo the boundary operators  $(b_1u, b_2u) = (u\vert_{\Sigma_+},u\vert_{\Sigma_-})$
whose domain   is
\begin{eqnarray*}D(b_1,b_2) = \lbrace u\in L^p(\Omega); (u\vert_{\Sigma_+},u\vert_{\Sigma_-})\in L^2(\Sigma_+ ; \vert \vec a \cdot\hat n \vert d\sigma)\times L^2(\Sigma_- ; \vert \vec a \cdot\hat n \vert d\sigma)\rbrace
\end{eqnarray*}
In this case we get the following 
 \begin{theorem} Assume $p\ge 2$, and let $f\in L^q$ ($\frac{1}{p}+\frac{1}{q}=1$) and $a_0\in L^{\infty}(\Omega)$.
Suppose there exists $\tau\in C^1(\bar\Omega)$ such that
\begin{equation}
\label{semi}
\vec{a}\cdot{\vec \nabla}\tau +\frac{1}{2}(\vec \nabla\cdot\vec{a} + a_0)\ge 0.
\end{equation}
Define
\[\phi(u) := \frac{1}{p}\int_\Omega e^{2\tau}\vert e^{-\tau}u\vert^pdx + \frac{1}{2}\int_\Omega\vec{a}\cdot{\vec \nabla}\tau\vert u\vert^2dx +\frac{1}{4}\int_\Omega(\vec \nabla\cdot\vec{a} + a_0)\vert u\vert^2dx +\int_\Omega ufdx,  \]
and
\[I(u) := \phi(u) + \phi^*(Au) + \frac{1}{2}\int_{\Sigma_+}{\vert u\vert^2\vert \vec a\cdot\hat n\vert d\sigma} + \frac{1}{2}\int_{\Sigma_-}{\vert u\vert^2\vert \vec a\cdot\hat n\vert d\sigma}\]
on the set
\[S :=\lbrace u\in L^p;\, u\in D(A)\cap D(b_1, b_2)\rbrace\]
\begin{enumerate}
\item Then there exists $\bar u \in S$ such that  $ I(\bar u) = \inf\lbrace I(u)\mid u \in S\rbrace = 0.$
\item The function $\bar v := e^{-\tau}\bar u$ satisfies the nonlinear transport equation

\begin{eqnarray}
\label{transport.plus}
\vec{a}\cdot\vec \nabla{\bar v}-\frac{a_0}{2}\bar v &=&\bar v\vert{\bar v}\vert^{p-2}+f \quad \mbox{\rm on }\Omega\\
\bar v &=&0 \quad \quad \quad \quad \quad \mbox{\rm on }\Sigma_+\nonumber
\end{eqnarray}

\end{enumerate}
\end{theorem}
\noindent{\bf Proof:} Consider
\begin{equation}
\label{LA}
M(u,p)=
\begin{cases}
{\phi(u) + \phi^*(Au +p) + \ell(b_1(u), b_2(u))\quad \mbox{ if }u\in S}\cr {+\infty\quad \quad \quad \quad \quad \quad \quad \quad \quad \mbox{ if }u \notin S}\cr
\end{cases}
\end{equation}
Where $\ell(\cdot,\cdot): L^2(\Sigma_+ ; \vert \vec a\cdot\hat n \vert d\sigma)\times L^2(\Sigma_- ; \vert \vec a \cdot\hat n \vert d\sigma)\to \R$ is defined by
\[\ell(h,k) := \frac{1}{2}\int_{\Sigma_+}{\vert h\vert^2\vert \vec a\cdot\hat n\vert d\sigma} + \frac{1}{2}\int_{\Sigma_-}{\vert k\vert^2\vert \vec a\cdot\hat n\vert d\sigma}.\]
Since $A$ is skew-adjoint modulo the boundary, we conclude from Proposition 2.1  that $M$ is an ASD Langrangian on the space $L^p(\Omega)\times L^q(\Omega)$ and satisfies all the hypothesis of Proposition 3.1. There exists then $\bar u\in L^p(\Omega)$ such that
$0 = M(\bar u,0) = \inf\lbrace M(u,0); u\in L^p\rbrace$ 
which means that
$0 = I(\bar u) = \inf\lbrace I(u)\mid u\in S\rbrace$ and 
assertion 1) is verified.\\
To get 2) we observe again that
\begin{eqnarray*}
0 &=&I(\bar u)=\phi(\bar u) + \phi^*(A\bar u) + \frac{1}{2}\int_{\Sigma_+}{\vert \bar u\vert^2\vert \vec a\cdot\hat n\vert d\sigma} + \frac{1}{2}\int_{\Sigma_-}{\vert \bar u\vert^2\vert \vec a\cdot\hat n\vert d\sigma}\\
&\geq & \langle \bar u,A\bar u\rangle + \frac{1}{2}\int_{\Sigma_+}{\vert \bar u\vert^2\vert \vec a\cdot\hat n\vert d\sigma} + \frac{1}{2}\int_{\Sigma_-}{\vert \bar u\vert^2\vert \vec a\cdot\hat n\vert d\sigma}\\
&&= \int_{\Sigma_+}{\vert \bar u\vert^2\vert \vec a\cdot\hat n\vert d\sigma}\geq 0
\end{eqnarray*}
 In particular, $\phi(\bar u) + \phi^*(A\bar u)=\langle u,Au\rangle$ and $\int_{\Sigma_+}{\vert \bar u\vert^2\vert \vec a\cdot\hat n\vert d\sigma}=0$, in such a way that $A\bar u\in\partial\phi(\bar u)$ and $u\vert_{\Sigma_+} = 0$. 
In other words, 
\[\vec a\cdot\vec \nabla\bar u +\frac{1}{2}(\vec \nabla\cdot\vec a)\bar u = e^\tau\bar u\vert e^{-\tau} \bar u\vert^{p-2} + (\vec a\cdot\vec \nabla\tau) \bar u + \frac{1}{2}(\vec \nabla\cdot\vec a + a_0)\bar u +f\]
and $\bar u\vert_{\Sigma_+} = 0$. 
Multiply now both equations by $e^{-\tau}$ and use the product rule for differentiation to get
$\vec{a}\cdot\vec \nabla{\bar v}-\frac{a_0}{2}\bar v =\bar v\vert{\bar v}\vert^{p-2}+f$ and $\bar v\vert_{\Sigma_+} \equiv 0$, where $\bar v := e^{-\tau}\bar u$.\\

\noindent {\bf Case 2:  $1<p \leq 2$}\\
In this case, the right space is $X = L^2(\Omega)$ and $A:D(A) \to X^*$ is defined as in the case when $p\geq2$ but this time on the domain
\[D(A)=\lbrace u\in L^2(\Omega)\mid \vec a\cdot\vec \nabla u \in L^2(\Omega)\rbrace.\]
By Lemma 2.4, A is again skew-adjoint modulo the boundary operators  
$(b_1u, b_2u) = (u\vert_{\Sigma_+},u\vert_{\Sigma_-})$
whose domain   is
\[
D(b_1,b_2) = \lbrace u\in L^2(\Omega) \vert (u\vert_{\Sigma_+},u\vert_{\Sigma_-})\in L^2(\Sigma_+ ; \vert \vec a\cdot\hat n \vert d\sigma)\times L^2(\Sigma_- ; \vert \vec a \cdot\hat n \vert d\sigma)\rbrace.
\]
We then  obtain the following theorem:
\begin{theorem} Assume $1<p\leq 2$ and let $f\in L^2$ and $a_0\in L^{\infty}(\Omega)$.
Suppose there exists $\tau\in C^1(\bar\Omega)$ such that for some $ \epsilon >0$ we have:
\begin{equation}
\label{semiplus}
\vec{a}\cdot{\vec \nabla}\tau +\frac{1}{2}(\vec \nabla\cdot\vec{a} + a_0)\geq \epsilon>0.
\end{equation}
Define
\[\phi(u) := \frac{1}{p}\int_\Omega e^{2\tau}\vert e^{-\tau}u\vert^pdx + \frac{1}{2}\int_\Omega\vec{a}\cdot{\vec \nabla}\tau\vert u\vert^2dx +\frac{1}{4}\int_\Omega(\vec \nabla\cdot\vec{a} + a_0)\vert u\vert^2dx+\int_\Omega ufdx \]
and
\[I(u) := \phi(u) + \phi^*(Au) + \frac{1}{2}\int_{\Sigma_+}{\vert u\vert^2\vert \vec a\cdot\hat n\vert d\sigma} + \frac{1}{2}\int_{\Sigma_-}{\vert u\vert^2\vert \vec a\cdot\hat n\vert d\sigma}\]
On the set
\[S :=\lbrace u\in L^2;\,  u\in D(A)\cap D(b_1, b_2)\rbrace\]
Then
\begin{enumerate}
\item  There exists then $\bar u \in S$ such that $ I(\bar u) = \inf\lbrace I(u)\mid u \in S\rbrace = 0$
\item The function $\bar v := e^{-\tau}\bar u$ satisfies the nonlinear transport equation (\ref{transport.plus})
 \end{enumerate}
\end{theorem}
\noindent{\bf Proof}:
Define $M$ as in (\ref{LA})
 and again by Lemma 2.4, $A$ is skew-adjoint modulo the boundary. Since now $\phi $ is bounded on the bounded sets of $L^2$, we can now invoke Proposition 2.1 to conclude that $M(u,p)$ is an ASD Langrangian on the space $L^2(\Omega)\times L^2(\Omega)$. But in this case, 
$\phi$ is coercive because of Condition (\ref{semiplus}),
and therefore $\phi^*$ is bounded on bounded sets. All the hypothesis of Proposition 3.1 are now satisfied so there exists $\bar u\in L^2(\Omega)$ such that $0 = M(\bar u,0) = \inf\lbrace M(u,0)\mid u\in L^2\rbrace$. The rest follows as in the case when $p\geq 2$. \\

\noindent {\bf 2. Transport equation with a diffusion term}\\
In this case, the conditions on the smooth vector field $\vec a $ and $\Omega$ need not be as restrictive as in the case where the equation is purely governed by the transport operator. This is because the setting will only require that the operator $A : D(A) \to X^*$ defined by $u\mapsto \vec a\cdot\vec \nabla u + \frac{1}{2}(\vec \nabla\cdot\vec a)u$ 
be only anti-symmetric. The setup is as follows:\\
Let $X=L^2(\Omega)$  and consider the above operator $A$  on the domain
\[D(A) = \lbrace u \in L^2(\Omega)\mid \vec a\cdot\vec \nabla u + \frac{1}{2}(\vec \nabla\cdot\vec a)u \in L^2(\Omega)\rbrace.\]
 We then have the following 
\begin{theorem}
Assume $p\geq2$ and $m > 1$ and let $\vec a \in C^{\infty}(\bar\Omega)$ be a smooth vector field. Suppose  $a_0 \in L^{\infty}$ satisfies the following coercivity condition:
\begin{equation}
\label{coercive}
\frac{1}{2}(\vec \nabla\cdot\vec a + a_0) \geq 0.
\end{equation}
Consider the following convex and lower semi-continuous functional on $L^2(\Omega)$:
\[\phi(u) :=\begin{cases}
{\frac{1}{p}\int_\Omega\vert\nabla u\vert^pdx +\frac{1}{4}\int_\Omega(\vec \nabla\cdot\vec a + a_0)\vert u\vert^2dx +\frac{1}{m}\int_\Omega \vert u\vert^mdx +\int_\Omega ufdx \quad \mbox{\rm  if }u\in W^{1,p}_0(\Omega)}\cr {+\infty\quad \quad \quad \quad \quad \quad \quad \quad \quad \quad \mbox{ if }u \notin W^{1,p}_0(\Omega)}\cr
\end{cases}\]
 and define the functional $I(u) := \phi(u) +\phi^*(Au)$ on $W^{1,p}_0(\Omega)$. Then
\begin{enumerate}
\item  There exists $\bar u\in W^{1,p}_0(\Omega)$ such that $I(\bar u)=\inf\lbrace I(u); \,  u\in W^{1,p}_0(\Omega)\rbrace =0$. 
\item The minimizer $\bar u$ satisfies the equation
\begin{eqnarray*}
-\Delta_p \bar u + \frac{1}{2}a_0\bar u + \bar u\vert u\vert^{m-2} +f &=& \vec a\cdot\vec \nabla\bar u \quad \quad   {\rm on}\, \Omega\\
\bar u&=&0\quad \quad  {\rm on}\quad  \partial\Omega
\end{eqnarray*}
\end{enumerate}
\end{theorem}
\noindent {\bf Proof:} The functional $\phi $ has a symmetric domain that is contained in the domain of $A$. So now  Lemma 2.5 applies and  the Lagragian
\[M(u,p) := \begin{cases}{\phi(u) + \phi^*(Au -p) \mbox{ if }u\in W^{1,p}_0(\Omega)}\cr {+\infty\quad \quad \quad \mbox{ if }u \notin W^{1,p}_0(\Omega)}\cr
\end{cases}\]
is ASD. Now $\phi$ is obviously coercive on $L^2(\Omega)$  and therefore $\phi^*$ is bounded on bounded sets of $L^2(\Omega)$. Proposition 3.1 applies and we find  $\bar u\in L^2(\Omega)$ such that $0=M(\bar u,0)=\inf\lbrace M(u,0)\mid u\in L^2(\Omega)\rbrace$. 
Clearly  $\bar u\in W^{1,p}_0(\Omega)$ and the rest follows as in the preceeding cases.   

\subsection{ASD Lagrangians and maximal monotone operators}
 
Assuming that we are in a Hilbert space setting, namely $X = H = X^*$, and recalling that an ASD Lagrangian $L$ satisfies $L(x,p)\geq \langle x,p\rangle$ for every $(x,p)\in H\times H$, we can  
 consider the problem  of minimizing for each fixed $p$, the functional $I_p(x)=L(x,p) +\langle x,p\rangle $ over $x\in H$. The same proof as Proposition 3.1 above (established in \cite{G2}) yields  that if $p\mapsto L(0,p)$ is bounded on  every ball of $H$, then for each $p\in H$ the minimization problem
\[ 
 \inf\limits_{x\in H}I_p(x)=\inf\limits_{x\in H}\{L(x,p) +\braket{x}{p}\}
 \]
is equal to zero and is attained at some point $\bar x (p)\in H$. 
  
  If now $L$ is strictly convex in the first variable, then $\bar x (p)$ is unique, and therefore  we can define a map $X(\cdot): H\to H$ by making $X(p)$ the unique point satisfying $L(X(p),p) +\braket{X(p)}{ p} = 0$. The convexity of $L$ allows us to using a monotonicity argument and then show that the map $X$ is monotone, that is
  \[
 \hbox{ $\langle X(p)- X(q), p-q\rangle \geq 0$ for all $(p,q)\in H\times H$.}
  \]
  Furthermore, if $L(\cdot,\cdot)$ is uniformly convex in the first variable, that is if 
$L(x,p)-\frac{\var\| x\|^2}{2}$  is convex in $x$,  then one can show that $p\to X(p)$ is a Lipchitz continuous operator (See next section).  A standard argument using contraction mapping theorem then shows that under these assumptions, $X$ is then maximal monotone \cite{Br}. The following proposition summarises this discussion. The proof is a straightforward application of standard convex analysis results some of which discussed in the next section and in  \cite{G2}. The details are left to the interested reader.

\begin{proposition}
Let $L: H\times H\to\R$ be an anti-selfdual Lagrangian on a Hilbert space $H$ that is strictly convex in the first variable. Suppose that the map $p\mapsto L(0,p)$ is bounded on the bounded sets of $H$. Then there exists a monotone map $X: H\to H$ such that  
\[
0=L(X(p),p) +\braket{X(p)}{p}=\inf\limits_{x\in H}\{L(x,p)+\braket{x}{p}\}.
\]
Furthermore, if $L$ is uniformly convex in the first variable, then $X$ is a Lipschitz map which is  maximal monotone. 
\end{proposition}
\section{A variational  principle for general evolution equations}
In this section we develop further the variational theory for dissipative evolution equations via the theory of ASD Lagrangians. The goal is to extend the variational theory of gradient flows  \cite{GT1} and other parabolic equations  developed in \cite{G2} so as to include evolutions of the form 
 \begin{equation}
\label{main.evolution}
-\dot x(t) \in \partial \phi(x(t)) + Ax(t) + \omega x(t), 
\end{equation}
where $A$ is an unbounded positive operator and $\omega$ is any real number. The starting point is Theorem 4.2  of \cite{G2} (see also \cite{GT1}), which based on the fact that if $L:[0,T]\times H\times H\to \R$ is a (time-dependent) anti-self dual Lagrangian on a Hilbert space $H$, then it ``lifts" to an anti-selfdual Lagrangian on path space 
 \[
A^{2}_{H} = \{ u:[0,T] \rightarrow H; \,\dot{u} \in L^{2}_{H}  \}
\]
consisting of all absolutely continuous arcs $u: [0,T]\to H$, 
equipped with the norm
\[
     \|u\|_{A^{^{2}}_{H}} = (\|u(0)\|_{H}^{2} +
     \int_0^T \|\dot{u}\|^{2} dt)^{\frac{1}{2}}.
\]

By associating an appropriate Lagrangian $L$ to the convex functional $\phi$ and the operator $A$, as in the preceeding sections, then one can already deduce from that theorem,  a variational resolution for (\ref{main.evolution}), at least for $\omega =0$. However, the boundedness condition for Theorem 4.2  of \cite{G2} is too stringent for most applications, but can be considerably relaxed  when the Lagrangian $L(x, p)$ is autonomous.  The rest of the section consists of doing just that through a Yosida-type  $\la$-regularization procedure reminiscent of the standard theory for convex functions, which seems to apply as naturally to ASD Lagrangians. We shall also deal with the case when $\omega$ is not zero, since it will allow --among other things-- to relax the convexity assumptions on $\phi$. \\

Given a Hilbert space $H$, we recall that $L:[0,T]\times H\times H\to \R$ is said to be a (time-dependent) anti-self dual Lagrangian if for each $t\in [0,T]$, the Lagrangian $L(t,\cdot, \cdot)\in \ASDH$. On the other hand, a {\it self-adjoint boundary Lagrangian} $\ell$ is a convex lower semi-continuous functional on $H\times H$ that satisfies $\ell^*(x,p)=\ell(-x,p)$. 
 
Define the {\it Partial Domain} of $\partial L$ to be the set:
\[
{\rm Dom}_{1}(\partial L)=\{x\in X;\hbox{\rm  there exists $p,q\in X^{*}$ such that $(p,0)\in \partial L(x, q)$\}.}
\]
Note that if $L(x,p)=\phi (x) +\phi^{*}(-p)$ with $0$ assumed to be in the domain of $\partial \phi$, then $x_0$ belongs to  ${\rm Dom}_{1}(\partial L)$ if and only if it belongs to the domain of $\partial \phi$. 
 
We shall say that a Lagrangian $L:H\times H\to\overline{\R}$  is {\it uniformly convex in the first variable (resp. second variable)}  if there exists $\var >0$ such that for all $p\in H$ (resp. for all $x\in H$) the Lagrangian 
\[
 L(x,p)-\frac{\var\| x\|^2}{2} \left(\mbox{resp. }L(x,p)-\frac{\var\| p\|^2}{2}\right)
 \]
 is convex in $x$ (resp. in $p$).  Here is the main result of this section.

\begin{theorem} Let $L$ be an autonomous anti-selfdual Lagrangian on a Hilbert space $H\times H$ that is uniformly convex in the first variable. Assuming ${\rm Dom}_{1}(\partial L)$ is non-empty, then for any $\omega \in {\bf R}$,  there exists a semi-group of operators $(T_{t})_{t\in {\bf R}^{+}}$ on $H$   such that $T_{0}=Id$ and for any $x_0\in {\rm Dom}_{1}(\partial L)$, the path $x(t)=T_tx_0$ satisfies the following:
\begin{equation}
\label{eqn:ultimate}
- ({\dot x}(t)+\omega x(t), x(t))\in \partial L(x(t), 
{\dot x} (t)+\omega x(t))
\end{equation}
  The path $(T_tx)_t$ is obtained as the unique minimizer on $A_{H}^2$ of the functional
\[
 \tilde I(u)=  \int_0^T e^{2\omega t} L( u(t), e^{-\omega t}(\frac{d}{dt}e^{\omega t} {u}(t)))dt +\frac{1}{2}\|u(0)\|^{2} -2\langle x, u(0)\rangle +\|x_0\|^{2} +\frac{1}{2}\|e^{\omega T}u(T)\|^{2}.
\]
Moreover, $\tilde I(x)=\inf\limits_{u\in A_{H}^2}\tilde I(u)=0$.\\ 
The semi-group is:
\begin{enumerate}
\item a contraction when $\omega > 0$
\item 1-Lipschitz when $\omega = 0$
\item locally Lipschitz in $t$ when $\omega < 0$
\end{enumerate} 
\end{theorem} 
 
First we shall  prove the following improvement of Theorem 4.2 of \cite{G2} provided $L$ is  autonomous. The boundedness condition is still there, but we first cover the semi-convex case.
 \begin{proposition}
Assume  $L:H\times H\to\overline{\R}$  is an autonomous anti-selfdual Lagrangian that is uniformly convex and  suppose 
\begin{equation}
L(x,0)\le C(\| x\|^2+1) \quad \hbox{\rm for all $x\in H$. }
\end{equation}
Then, for any $w\in\R$ and any $x_0\in H$, there  exists $\hat x \in C^1\big( [0,T]:H\big)$ such that $\hat x(0)=x_0$ and 
\begin{equation}
\label{P1}
\int_0^T e^{2\omega t} L(e^{-\omega t} {\hat x}(t), e^{-\omega t} {\dot {\hat x}}(t))dt +\frac{1}{2}\|{\hat x}(0)\|^{2} -2\langle x_0, {\hat x}(0)\rangle +\|x_0\|^{2} +\frac{1}{2}\|{\hat x}(T)\|^{2}=0
\end{equation}
\begin{equation}
\label{P2}
-e^{-wt}\big( \dot{\hat x}(t),\hat x(t)\big)\in \partial L\big( e^{-wt}\hat x(t),e^{-wt}\dot{\hat x}(t)\big)
       \quad\forall t\in [0,T], 
\end{equation}
\begin{equation}
\label{P3}
\|\dot{\hat x}(t)\|\le C(w,T)\|\dot{\hat x}(0)\|   \quad\forall t\in [0,T], 
\end{equation}
where $C(w,T)$ is a positive constant.
 \end{proposition}
We shall need first the following lemmas which show how uniform convexity of the Lagrangian yield certain regularity properties of the solutions.  
  
\begin{lemma}
Let $F:H\to \overline{\R}$ be convex and lower semi-continuous and such that its Legendre dual  $F^*$ is uniformly convex. Then for every $x\in H$, the subdifferential $\partial F (x)$ is nonempty, is single-valued and the map $x\to \partial F (x)$ is Lipschitz on $H$. 
 \end{lemma}

\noindent {\bf Proof:} Since $F^*$ is uniformly convex, then  $F^*(x)=G(x)+\frac{{\var \| x\|}^2}{2}$ for some   convex  lower semi-continuous function $G$ and some $\epsilon>0$. It follows that 
 $F^*(x)\ge C+\braket{a}{x}+\frac{{\var\| x\|}^2}{2}$ for some $a\in H$ and $C>0$, hence 
$F(x)=F^{**}(x)\le C(1+\| x\|^2)$ which means that $F$ is 
 sub-differentiable for all $x\in H$. 

Consider now $p_j\in  \partial F(x_j)$ for $ j=1,2$ in such a way that 
$x_j \in \partial F^*(p_j)=\partial G(p_j)+\var p_j$.
 By monotonicity, we have
$0 \le \braket{p_1-p_2}{\partial G(p_1)-\partial G(p_2)}
=\braket{p_1-p_2}{x_1-\var p_1-x_2+\var p_2}
$
which yields that $ \var\| p_1-p_2\| \le \| x_1-x_2\|$ and we are done.
 
\begin{lemma} Assume  $L:H\times H\to\overline{\R}$  is an anti-selfdual Lagrangian that is uniformly convex in both variables. Then, for all $x,u\in H$, there  exists  a unique $v\in X$ --denoted $v=R(u,x)$ such that $x= \partial_2L(u,v)$. Moreover,  the map  $(u,x)\to R(u,x)$ is jointly Lipschitz on $H\times H$.
 \end{lemma}
\noindent {\bf Proof:}   Since $L$ is uniformly convex, then 
$
L(x,p)=M(x,p)+\var\left(\frac{\| x\|^2}{2}+\frac{\| p\|^2}{2}\right)
$, 
where $M$ is convex lower semi-continuous, in such a way that  
 $x=\partial_2L(u,v)$ if and only if $0\in \partial_2M(u,v)+\var v-x$ if and only if 
 $v$ is the solution to the following minimization problem
\[
\min_p\left\{ M(u,p)+\frac{\var\| p\|^2}{2}-\braket{x}{p}\right\}. 
\]
But for each fixed $u$ and $x$, the map $p\mapsto M(u,p)-\braket{x}{p}$
majorizes a linear functional  and therefore the minimum is attained uniquely at $v$ by strict convexity and obviously  $x=\partial_2L(u,v)$.\\
To establish the Lipschitz property, write 
\[
R(u_1,x_1)-R(u_2,x_2)=R(u_1,x_1)-R(u_1,x_2)+R(u_1,x_2)-R(u_2,x_2).
\]
We first  bound $\| R(u_1,x_1)-R(u_1,x_2)\|$ as follows:

Since $x_1 =\partial_2L\big( u_1,R(u_1,x_1)\big)$, $x_2 =\partial_2L\big( u_1,R(u_1,x_2)\big)$ and 
$ L(u_1,v)=M(u_1,v)+\frac{\var\| v\|^2}{2}$ for some $M$ convex and lower semi-continuous, it follows that 
 $x_j =\partial_2M\big( u_1,R(u_1,x_j)\big) +\var R(u_1,x_j)$ for $j= 1,2$,  so by monotonicity we get
\begin{eqnarray*}
0  &\le & \left\langle R(u_1,x_1)-R(u_1,x_2),\partial M\big( u_1,R(u_1,x_1)\big) 
    -\partial M\big( u_1,R(u_1,x_2)\big)\right\rangle\\
&=& \left\langle R(u_1,x_1)-R(u_1,x_2),x_1
   -\var R(u_1,x_1)-x_2+\var R(u_1,x_2)\right\rangle
\end{eqnarray*}
which yields that 
\[
\var {\left\| R(u_1,x_1)-R(u_1,x_2)\right\|}^2
    \le {\left\| R(u_1,x_1)-R(u_1,x_2)\right\|} 
    \| x_1-x_2\| 
\]
and therefore
\begin{equation}
\label{bound.1}
\left\| R(u_1,x_1)-R(u_1,x_2)\right\|
    \le\frac{1}{\var} \| x_1-x_2\|.
\end{equation}
Now we bound $\left\| R(u_1,x_2)-R(u_2,x_2)\right\|$.

Let $x_2 = \partial_2L\big( u_j,R(u_j,x_2)\big)=\partial_2M\big( u_j, R(u_j,x_2)\big) +\var R(u_j,x_2)$ for $ j=1,2$.
 and write by monotonicity that
\begin{eqnarray*}
0\le \left\langle R(u_1,x_2)-R(u_2,x_2),\partial_2M\big( u_1,R(u_1,x_2)\big) 
   -\partial_2M\big( u_1,R(u_2,x_2)\big)\right\rangle
\end{eqnarray*}
Setting $ p_j=R(u_j,x_2)$, we have with this notation
 
\begin{eqnarray*}
\langle p_1-p_2,\partial_2M(u_1,p_2)-\partial_2M(u_2,p_2)\rangle &\le&
\langle p_1-p_2,\partial_2M(u_1,p_1)-\partial_2M(u_2,p_2)\rangle \\
&=&\langle p_1-p_2,x_2-\var p_1-x_2+\var p_2\rangle \\
&=&-\var \| p_1-p_2\|^2.
\end{eqnarray*}
so that
$\var \| p_1-p_2\|^2\le \| p_1-p_2\| \|\partial_2M(u_1, p_2)-\partial_2M(u_2,p_2)\|$, and since 
$\partial_2M(u_j,p_2)=\partial_2L(u_j,p_2)-\var p_2$, we get that
\[
 \var \| p_1-p_2\| 
\le \| \partial_2L(u_1,p_2)-\partial_2L(u_2,p_2)\|
\le \| \partial L(u_1,p_2)-\partial L(u_2,p_2)\| 
\]
Here we use the fact that $L$ is both anti-selfdual  and uniformly convex, to deduce that  $L^*$ is also uniformly convex. We then apply Lemma 4.2 to get:
\[
 \| \partial L(u,p)-\partial L(u',p')\| \le C\big( \| u-u'\| +\| p-p'\|\big)
 \]
 from which follows   that
$ \| p_1-p_2\|\le \frac{C}{\epsilon}\| u_1-u_2\|$, hence 
\begin{equation}
\label{bound.2}
 \| R(u_1,x_2)-R(u_2,x_2)\|\le  \frac{C}{\epsilon}\| u_1-u_2\|.
\end{equation}
Combining estimates (\ref{bound.1}) and (\ref{bound.2}), we finally get
 \[
 \| R(u_1,x_1)-R(u_2,x_2)\|\le
\frac{1}{\epsilon}(1+C)\left( \| u_1-u_2\| +\| x_1-x_2\|\right).
\]

We can now deduce the following corollary that gives a regularity result for certain flows.
\begin{lemma}
Assume  $L:H\times H\to\overline{\R}$  is an anti-selfdual Lagrangian that is uniformly convex in both variables. Suppose the paths $v,x,u:[0,T]\to H$ are such that  $x,u\in C([0,T]; H)$ and $-x(t)=\partial_2L\big( u(t),v(t)\big)$ for almost all $t\in [0,T]$. Then $v \in C([0,T]; H)$ and $-x(t)=\partial_2L\big( u(t),v(t)\big)$ for all $t\in [0,T]$.
\end{lemma}

\noindent {\bf Proof of Proposition 4.1:}  Apply Theorem 4.2 of \cite{G2} to the Lagrangian $M(t,x,p)=e^{2\omega t} L(e^{-\omega t} {x}, e^{-\omega t}  p)$ which is also anti-selfdual (See \cite{G2}). There exists then 
 $\hat x \in A_{H}^2$ such that
$\big(\hat x (t),\dot{\hat x} (t)\big)\in \mbox{\rm Dom} (M)$ for almost all 
$t\in [0,T]$    and $ I( {\hat x})=\inf\limits_{u\in A_{H}^2}I(u)=0$, 
where 
\[
 I (u)=  \int_0^T M(t, u(t),\dot{u}(t))dt +\frac{1}{2}\|u(0)\|^{2} -2\langle x_0, u(0)\rangle +\|x_0\|^{2} +\frac{1}{2}\|u(T)\|^{2}.
\]
The  path ${\hat x}$ then satisfies:
${\hat x}(0)=x_{0}$ and for almost all $t\in [0,T]$,
\[
-\big(\dot{\hat x}(t),\hat x(t)\big)\in \partial M\big( t,\hat x(t),\dot{\hat x}(t)\big) 
\]
and the chain rule  
\[
\partial M(t,x,p)= e^{wt}\partial L\big( e^{-wt}x,e^{-wt}p\big)
\]
to get that for almost all $t\in [0,T]$
\[
 -e^{-wt}\big(\dot{\hat x},\hat x(t)\big)\in \partial L
\big(e^{-wt}\hat x(t),e^{-wt}\dot{\hat x}(t)\
 \]
Apply Lemma 4.4 to $x(t)=u(t)=e^{-wt}\hat x(t)$ and $v(t)=e^{-wt}\dot{\hat x}(t)$ to conclude that  $\dot{\hat x}\in C\big( [0,T]:H\big)$. Thus $\hat x\in C^1\big( [0,T]:H\big)$. Since $L$  is anti-selfdual and uniformly convex, we get from Lemma 4.2  that $ (x,p)\mapsto \partial L(x,p)$  is Lipschitz. So by continuity, we have now for all $t\in [0,T]$ 
\[ 
-e^{-wt}\big(\dot{\hat x}(t),\hat x(t)\big)\in \partial L\big(e^{-wt}\hat x(t),e^{-wt}\dot{\hat x}(t)\big)
\]
and (\ref{P2}) is verified.

To establish (\ref{P3}), we first differentiate to obtain:
\[
e^{-2wt}\frac{d}{dt} \|\hat x(t)-e^{-wh}\hat x(t+h)\|^2=
2e^{-2wt}\braket{\hat x(t)-e^{-wh}\hat x(t+h)}{\dot{\hat x}(t)-e^{-wh}\dot{\hat x}(t+h)}.
\]
Setting now
$v_1(t) =\partial_1L\big( e^{-wt}\hat x(t),e^{-wt}\dot{\hat x}(t)\big)$ and 
$v_2(t) =\partial_2L\big(e^{-wt}\hat x(t),e^{-wt}\dot{\hat x}(t)\big)$, we obtain 
 from (\ref{P2})  and monotonicity that 
\begin{eqnarray*}
e^{-2wt}\frac{d}{dt} \|\hat x(t)-e^{-wh}\hat x(t+h)\|^2&=&
 \langle e^{-wt}\hat x(t)-e^{-w(t+h)}\hat x(t+h),-v_1(t)+v_1(t+h)\rangle +\\
& &\quad \langle e^{-wt}\dot{\hat x}(t)-e^{-w(t+h)}\dot{\hat x}(t+h),-v_2(t)+v_2(t+h)\rangle\\
 &\leq& 0.
\end{eqnarray*}
We conclude from this that
\begin{eqnarray*}
 \frac{\|\hat x(t)-e^{-hw}\hat x(t+h)\|}{h}\leq \frac{\|\hat x(0)-e^{-hw}\hat x(h)\|}{h}
\end{eqnarray*}
and as we take $h \to 0$, we get 
$
\| w\hat x(t)+\dot{\hat x}(t)\| \leq 
\| w x_0+\dot{\hat x}(0)\| + \|\omega x_0\|.$ 
 Therefore
 \[
 \|\dot{\hat x}(t)\|\le \|\dot{\hat x}(0)\| +|w|\|\hat x(t)\| + T\|\omega x_0\|
 \]
  and 
\[
   \|\hat x(t)\|\le\int_0^t\|\dot{\hat x}(s)\|\, ds\le\|\dot{\hat x}(0)\| T+|w|\int_0^t\|\hat x(s)\|\, ds.
   \]
   It follows from Gronwall's  inequality that 
$\|\hat x(t)\|\le\big(\|\dot{\hat x}(0)\|+ \|\omega x_0\|\big) \big( C+|w|e^{|w|T}\big)$ for all $t\in [0,T]$ 
and finally that $ \|\dot{\hat x}(t)\|\le\|\dot{\hat x}(0)\|\big(C+|w|+|w|^2e^{|w|T}\big)$. \\
 We now proceed with the proof of Theorem 4.1. For that we  associate a Yosida-type  $\la$-regularization of the Lagrangian so that the boundedness condition in Proposition 4.1 is satisfied, then we  make sure that all goes well when we take the limit as $\la$ goes to $0$. First, we need the following lemmas relating the properties of a Lagrangian to those of its $\la$-regularization.
 \begin{lemma}
For a convex functional $L:H\times H\to\overline{\R}$, define for  each $\la >0$, the Lagrangian
\begin{eqnarray*}
L_\la (x,p):=\inf_z\left\{ L(z,p)+\frac{\| x-z\|^2}{2\la}\right\} +\frac{\la\| p\|^2}{2}.
\end{eqnarray*}
\begin{enumerate}
\item If $L$ is anti-selfdual, then $L_\lambda$ is also anti-selfdual.
\item If $L$ is uniformly convex in the first variable, then $L_\lambda$  is uniformly convex (in both variables) on $H\times H$.

\end{enumerate}
\end{lemma}
 \noindent {\bf Proof:} Fix $(q, y)\in X^{*}\times X$ and write:
\begin{eqnarray*}
  (L_\lambda)^{*} (q,y)
  &=&\sup\{\langle q, x\rangle +  \langle y, p\rangle- L(z, p)- \frac{\| x-z\|^2}{2\la} -\frac{\la\| p\|^2}{2}; (z, x, p)\in X\times X\times X^{*}\}\\
&=&  \sup\{\langle q, v+z\rangle +  \langle y, p\rangle- L(z, p)- \frac{\| v\|^2}{2\la} -\frac{\la\| p\|^2}{2}; (z, v, p)\in X\times X\times X^{*}\}\\
&=&\sup_{p\in X^{*}} \left\{\langle y,p\rangle +\sup\limits_{(z, v)\in X\times X}\{\langle q,v+z\rangle -L(z,p)-\frac{\| v\|^2}{2\la}\}-\frac{\la\| p\|^2}{2}\right\}\\
&=&\sup_{p\in X^{*}}\left\{\langle y,p\rangle +\sup_{z\in X}\{\langle q,z\rangle -L(z,p)\} +\sup_{v\in X} \{\langle q,v\rangle -\frac{\| v\|^2}{2\la}\}-\frac{\la\| p\|^2}{2}\right\}\\
&=&\sup_{p\in X^{*}} \left\{\langle y,p\rangle +\sup_{z\in X}\{\langle q,z\rangle -L(z,p)\} +\frac{\la\| q\|^2}{2}-\frac{\la\| p\|^2}{2})\right\}\\
&=&\sup_{p\in X^{*}} \sup_{z\in X}\left\{\langle y,p\rangle +\langle q,z\rangle -L(z,p)-\frac{\la\| p\|^2}{2}\right\} +\frac{\la\| q\|^2}{2}\\
&=& (L+T)^{*}(q,y)+\frac{\la\| q\|^2}{2}
\end{eqnarray*}
where $T(z,p):= \frac{\la\| p\|^2}{2}$ for all $(z,p)\in X\times X^{*}$. Note now that 
\begin{eqnarray*}
T^*(q,y)=\sup_{z,p}\left\{\langle q,z\rangle +\langle y,p\rangle -  \frac{\la\| p\|^2}{2}\right\}=\left\{\begin{array}{lll}+\infty &\hbox{if }&q\ne 0\\ \frac{\|y\|^2}{2\la}&\hbox{if }&q=0\end{array} \right.
\end{eqnarray*}
in such a way that by using the duality between sums and convolutions in both variables, we get
\begin{eqnarray*}
(L+T)^{*}(q,y)&=&{\rm conv}(L^*,T^*)(q,y)\\
&=&\inf_{r\in X^{*},z\in X}\left\{ L^*(r,z)+T^*(-r+q,-z+y)\right\}\\
&=&\inf_{z\in X}\left\{ L^*(q,z)+\frac{\| y-z\|^2}{2\la}\right\}
\end{eqnarray*}
and  finally 
\begin{eqnarray*}
L_\la^{*} (q,y)&=&(L+T)^{*}(q,y)+ \frac{\la\| q\|^2}{2}\\
&=&\inf_{z\in X}\left\{ L^*(q,z)+\frac{\| y-z\|^2}{2\la}\right\}+ \frac{\la\| q\|^2}{2}\\
&=&\inf_{z\in X}\left\{ L(-z,-q)+\frac{\la\| q\|^2}{2}+\frac{\| y-z\|^2}{2\la}\right\}\\
  &=&L_\la(-y,-q).
\end{eqnarray*}

 (2) For each  $\la >0$, there exists $\var >0$ such that $ M(x,p):=L(x,p)-\frac{\var\| x\|^2}{\la^2}$ is convex. Pick $\delta =\frac{1-\frac{1}{1+\var}}{\la}$ so that $1+\var =\frac{1}{1-\la\delta}$ and write
 
\begin{eqnarray*}
L_\la (x,p)-\frac{\la\| p\|^2}{2}-\delta\frac{\| x\|^2}{2} &=&
\inf_z\left\{ L(z,p)+\frac{\| x-z\|^2}{2\la}-\frac{\delta\| x\|^2}{2}\right\}\\
&=&\inf_z\left\{ L(z,p)+\frac{\| x\|^2}{2\la}-\frac{\braket{x}{z}}{\la}+
   \frac{\| z\|^2}{2\la}-\frac{\delta\| x\|^2}{2}\right\}\\
&=&\inf_z\left\{ L(x,p)+\frac{{\left\| \sqrt{\frac{1}{\la}-\delta}x\right\|}^2}{2}
    -\frac{\braket{x}{z}}{\la}+\frac{\| z\|^2}{2\la}\right\}\\
&=&\inf_z\left\{ L(z,p)+\frac{{\left\|\sqrt{1-\la\delta}x\right\|}^2}{2\la}
    -\frac{\braket{\sqrt{1-\la\delta}x}{\frac{z}{\sqrt{1-\la\delta}}}}{\la}
    +\frac{\| z\|^2}{2\la}\right\}\\
&=&\inf_z\left\{ M(z,p)+\frac{\var\| z\|^2}{2\la}
   +\frac{{\left\|\sqrt{1-\la\delta}x\right\|}^2}{2\la}
   -\frac{\braket{\sqrt{1-\la\delta}x}{\frac{z}{\sqrt{1-\la\delta}}}}{\la}
    +\frac{\| z\|^2}{2\la}\right\}\\
&=&\inf_z\left\{ M(z,p)+\frac{(1+\var )\| z\|^2}{2\la}
   -\frac{\braket{\sqrt{1-\la\delta} x}{\frac{z}{\sqrt{1-\la \delta}}}}{\la}
   +\frac{{\left\|\sqrt{1-\la \delta}x\right\|}^2}{2\la}\right\}\\
&=&\inf_z\left\{ M(z,p)+\frac{{\left\|\frac{z}{\sqrt{1-\la\delta}}-\sqrt{1-\la\delta}x\right\|}^2}{2\la}\right\}
\end{eqnarray*}
 which means that 
$(z,p,x)\mapsto M(z,p)+\frac{{\left\| \frac{z}{\sqrt{1-\la\delta}}-\sqrt{1-\la\delta}x\right\|}^2}{2\la}$
 is convex and therefore the infimum in $z$
is  convex, which means that   
$L_\la (x,p)-\frac{\la\| p\|^2}{2}-\delta\frac{\| x\|^2}{2}
$
is itself convex, meaning that  $L_\la$ is uniformly convex.
 
\begin{lemma} For a given convex functional $L:H\times H\to\overline{\R}$ and  $\la >0$, denote for each $(p,x)\in H\times H$, by $J_\la (x,p)$ the minimizer of the following optimization problem:
\[
 \inf_z\left\{ L(z,p)+\frac{\| x-z\|^2}{2\la}\right\}.
 \]
 \begin{enumerate}
\item For each $(x,p)\in H\times H$, we have
 \begin{equation}
\label{Jlambda}
\partial_1L_\la (x,p)=\frac{x-J_\la (x,p)}{\la}\in \partial_1L\big( J_\la (x,p),p\big).
\end{equation}

\item If  $L:H\times H\to\overline{\R}$  is an anti-selfdual Lagrangian that is uniformly convex in the first variable, then the map $(x,p)\to J_\la (x,p)$ is Lipschitz on $H\times H$. 
\end{enumerate}
\end{lemma}

\noindent {\bf Proof:} (1) is straightforward. For (2), use Lemma 4.5 to deduce that $L_\la$ is  anti-selfdual and uniformly convex in both variables, which means that  $L_\la^*$ is also uniformly convex in both variables. It follows from Lemma 4.2  that $(x,p)\mapsto \partial L_\la (x,p)$ is Lipschitz. From (\ref{Jlambda}) above, we see that 
$J_\la (x,p)=x-\la \partial_1L_\la (x,p)$ is Lipschitz as well.\\

The following  lemma  will be useful in obtaining a uniform bound on the first derivatives of the family of approximate solutions. 

\begin{lemma}
Assume  $L:H\times H\to\overline{\R}$  is an anti-selfdual Lagrangian and let $L_\lambda$ be its $\lambda$-regularization, then the following hold:
\begin{enumerate}
\item If $ -(y,x)=\partial L_\la (x,y)$, then necessarily
$-\big( y,J_\la (x,y)\big)\in \partial L\big( J_\la (x,y),y\big).$
\item If $0\in {\rm Dom}_{1}(\partial L)$, then there exists a constant  $C>0$ such that $\| y_\la\|\le C$ whenever  $y_\la$ solves $-(y_\la ,0)=\partial L_\la (0,y_\la )$.
\end{enumerate}
\end{lemma}

\noindent {\bf Proof:}  (1) \quad If $ -(y,x)=\partial L_\la (x,y)$ then  $ L_\la (x,y)+L_\la ^*(-y,-x)=-2\langle x,y)$ and since 
$L$ is an ASD Lagrangian, we have $L_\la (x,y)+L_\la (x,y)=-2\langle x,y)$, hence 
 \begin{eqnarray*}
-2\braket{x}{y} &=& L_\la (x,y)+L_\la (x,y)\\
 &=&2\left( L\big( J_\la (x,y),y\big) 
+\frac{\| x-J_\la (x,y)\|^2}{2\la}+\frac{\la\| y\|^2}{2}\right)\\
&=& L^*\big( -y,-J_\la (x,y)\big) +L\big( J_\la (x,y),y\big) 
   +2\left(\frac{\|-x+J_\la (x,y)\|^2}{2\la}+\frac{\la\| y\|^2}{2}\right)\\
&\ge& -2\langle y,J_\la (x,y)\rangle +2 \langle -x+J_\la (x,y),y\rangle\\
&=& -2\langle x,y)
\end{eqnarray*}
The second last inequality is deduced by applying Fenchel's inequality to the first two terms and the last two terms. The above chain of inequality shows that all inequalities are equalities. This implies, again by Fenchel's inequality that $-\big( y,J_\la (x,y)\big)\in \partial L\big( J_\la (x,y),y\big)$.\\
  
(2) If $-(y_\la ,0)=\partial L_\la (0,y_\la )$, we get from Lemma 4.6.(1) that $-y_\la =\frac{-J_\la (0,y_\la )}{\la}\in \partial_1L\big( J_\la (0,y_\la ),y_\la\big)$, and by the first part of this lemma, that $-\big( y_\la ,J_\la (0,y_\la )\big)\in \partial L\big( J_\la (0,y_\la ),y_\la\big)$. 
 
 Now since $0\in {\rm Dom}_{1}(\partial L)$, there exists $\hat p$ such that $ \partial_1L(0,\hat p)\ne \emptyset$ and $0\in \partial_2L(0,\hat p)$.
Setting $ v_\la = J_\la (0,y_\la )$, and 
since  $-\big( y_\la ,v_\la )\big)\in \partial L\big( v_\la, y_\la\big)$, we get from monotonicity and by the fact that $y_\la=\frac{v_\la}{\la}$,
\begin{eqnarray*}
0 &\le&\braket{(0,\hat p)-(v_\la ,y_\la )}{\big( \partial_1L(0,\hat p),\partial_2L(0,\hat p)\big) -(-y_\la ,-v_\la )}\\
&=&\braket{(0,\hat p)-(v_\la ,y_\la )}{(\partial_1L(0,\hat p),0)-(\frac{-v_\la}{\la},-v_\la )}\\
  &=&-\frac{\| v_\la\|^2}{\la}-\braket{v_\la}{\partial_1L(0,\hat p)} 
    +\braket{\hat p}{v_\la} -\braket{y_\la}{v_\la}\\
  &=&-2\frac{\| v_\la\|^2}{\la}-\braket{v_\la}{\partial_1L(0,\hat p)}+\braket{\hat p}{v_\la}
\end{eqnarray*} 
which yields that $2\frac{\| v_\la\|}{\la}\le\| \partial_1L(0,\hat p)\| +\|\hat p\|$ and finally the desired bound
$2\| y_\la\|\le\| \partial_1L(0,\hat p)\| +\|\hat p\|$ for all $\la >0$. \\

\noindent {\bf End of Proof of Theorem 4.1:}  Let  $M_\la(t,x,p)=e^{2\omega t} L_\la(e^{-\omega t} {\hat x}, e^{-\omega t}  p)$  which is also anti-selfdual  and uniformly convex by Lemma 4.5. \\
We now have $L_\la (t, x,0)\le L(0,0)+\frac{\| x\|^2}{2\la}$, hence Proposition 4.2 applies and we get for all $\la >0$ a solution $ x_\la\in C^1\big([0,T]:H\big)$ such that $x_\la (0)=x_0$, 
\begin{eqnarray}
\label{P11}
\int_0^T M_\la \big( t,x_\la (t),\dot x_\la (t)\big)\, dt+\ell\big( x_\la (0),x_\la (T)\big) =0
\end{eqnarray}
\begin{eqnarray}
\label{P12}
\hbox{$-e^{-wt}\big(\dot x_\la (t),x_\la (t)\big)\in \partial L_\la \big(e^{-wt}x_\la (t),e^{-wt}\dot x_\la (t)\big)$ for all $t \in [0,T]$}
\end{eqnarray}
\begin{eqnarray}
\label{P13}
\|\dot x_\la (t)\|\le C(w,T)\|\dot x_\la (0)\|. 
\end{eqnarray}
Here $\ell\big( x_\la (0),x_\la (T)\big) =\frac{1}{2}\|x_\la(0)\|^{2} -2\langle x, x_\la(0)\rangle +\|x\|^{2} +\frac{1}{2}\|u_\la(T)\|^{2}.$
By the defintion of $M_\la(t,x,p)$, identity (\ref{P11}) can be written as
\begin{eqnarray}
\label{P15}
\int_0^T e^{2wt}L_\la\big( e^{-wt}x_\la (t),e^{-wt}\dot x_\la (t)\big)\, dt+\ell\big(x_\la (0),x_\la (T)\big) =0,
\end{eqnarray}
and since 
\begin{eqnarray*}
L_\la (x,p)=L\big( J_\la (x,p),p\big) +\frac{\| x-J_\la (x,p)\|^2}{2\la}+\frac{\la\| p\|^2}{2}
\end{eqnarray*}
Equation (\ref{P11}) can be written as
\begin{eqnarray*}
\int_0^T e^{2wt}\big( L(v_\la (t),e^{-wt}\dot x_\la (t)\big)
   +\frac{\| e^{-wt}x_\la (t)-v_\la (t)\|^2}{2\la}
   +\frac{\la\| e^{-wt}\dot x_\la (t)\|^2}{2})\, dt
   +\ell\big( x_\la (0),x_\la (T)\big)=0
\end{eqnarray*}
where $v_\la (t)=J_\la\big( e^{-wt}x_\la (t),e^{-wt}\dot x_\la (t)\big)$. Using Lemma 4.6.(1),  we get from  (\ref{P12}) that for all $t$,
\begin{eqnarray}
\label{P14}
-e^{-wt}\dot x_\la (t)=\partial_1L_\la \big(e^{-wt}x_\la (t),e^{-wt}\dot x_\la (t)\big)
=\frac{e^{-wt}x_\la (t)-v_\la (t)}{\la}
\end{eqnarray}
Setting $t=0$ in (\ref{P12})  we get $-\big(\dot x_\la (0),0\big)\in \partial L_\la \big( 0,\dot x_\la (0)\big)$, 
and since $0\in {\rm Dom}_1L$,  we can apply Lemma 4.7.2) to get that $\|\dot x_\la (0)\|\le C$ for all $\la >0.$ 
Now plug this inequality in  (\ref{P13}) and we obtain:
\[\|\dot x_\la (t)\|\le D(w,T)\quad\forall\la >0\quad\forall t\in [0,T]\]
This yields by (\ref{P14}) that
\[ \| e^{-wt}x_\la (t)-v_\la (t)\|\le e^{|w|T}D(w,T)\la\quad\forall t\in [0,T],\] 
hence
\begin{eqnarray}
\frac{\| e^{-wt}x_\la (t)-v_\la (t)\|^2}{\la}\to 0
\end{eqnarray}
uniformly in $t$. Moreover, since
${\|\dot x_\la (\cdot )\|}_{A_H^2}\le D(w,T)$ for all $\la >0$, 
there exists $\hat x \in A_H^2$ such that --up to a subsequence--
\begin{eqnarray}
x_\la \rightharpoonup\hat x \mbox{ in }A_H^2
\end{eqnarray}
and again by  (\ref{P14}) we have
\begin{eqnarray}
\int_0^T\|v_\la(t)-e^{-wt}\hat x(t)\|_H^2dt\to 0,
\end{eqnarray}
while clearly
\begin{eqnarray}
\la \frac{\| e^{-wt}\dot x_\la (t)\|^2}{2}\to 0\mbox{ uniformly}.
\end{eqnarray}
Now use (38)--(41) and the lower semi-continuity of $L$, to deduce from (\ref{P15}), that as $\la\to 0$ we have
\begin{eqnarray*}
I\big(\hat x \big) =
   \int_0^T e^{2wt}L\big( e^{-wt}\hat x(t), e^{-wt}\hat x(t)\big)\, 
   dt+\ell\big( \hat x (0),\hat x(T)\big)\le 0.
\end{eqnarray*}
Since we already know that 
$I \big( x \big)\ge 0$ for all $x \in A_H^2$, 
we finally get our claim that $0=I\big(\hat x \big)
   =\inf_{x \in A_H^2} I \big( x \big)$.\\
Now define $T_t x_0 := e^{-\omega t} \hat x(t)$.\\
It is easy to see that $x(t) := T_t x_0$ satisfies equation (23) and that $T_0 x_0 = x_0$. We need to check that $\lbrace T_t\rbrace_{t\in\R^+}$ is a semi-group. By uniqueness of minimizers, it is equivalent to show that for all $T' < T$, we have $w(t) := x(t+T')$  satisfies
\[
0 = \int_0^{T - T'} e^{2\omega t} L( w(t), e^{-\omega t}(\frac{d}{dt}e^{\omega t} {w}(t)))dt + \frac{1}{2}\|w(0)\|^{2} -2\langle T_{T'}x_0, w(0)\rangle +\|T_{T'}x_0\|^{2} +\frac{1}{2}\|e^{\omega (T-T')}w(T)\|^{2}
\]
By the defintion of $x(t)$ and the fact that $I\big(\hat x\big) = 0$ we have,
\[
0 = \int_0^T e^{2\omega t} L( x(t), e^{-\omega t}(\frac{d}{dt}e^{\omega t} {x}(t)))dt +\frac{1}{2}\|x(0)\|^{2} -2\langle x_0, x(0)\rangle +\|x_0\|^{2} +\frac{1}{2}\|e^{\omega T}x(T)\|^{2}
\]
Now let $0 < T' < T$. Since $x(t)$ satisfies equation (23) we have
\[
0 = \int_0^{T'} e^{2\omega t} L( x(t), e^{-\omega t}(\frac{d}{dt}e^{\omega t} {x}(t)))dt +\frac{1}{2}\|x(0)\|^{2} -2\langle x_0, x(0)\rangle +\|x_0\|^{2} +\frac{1}{2}\|e^{\omega T'}x(T')\|^{2}
\]
Subtract the two equations we get
\[
0 = \int_0^{T - T'} e^{2\omega t} L( x(t), e^{-\omega t}(\frac{d}{dt}e^{\omega t} {x}(t)))dt + \frac{1}{2}\|e^{\omega T} x(T)\|^2 - \frac{1}{2}\|e^{\omega T'} x(T')\|^2
\]
Make a substitution $s = t - T'$ and we obtain
\[
0 = \int_0^{T - T'} e^{2\omega t} L( w(t), e^{-\omega t}(\frac{d}{dt}e^{\omega t} {w}(t)))dt + \frac{1}{2}\|w(0)\|^{2} -2\langle T_{T'}x_0, w(0)\rangle +\|T_{T'}x_0\|^{2} +\frac{1}{2}\|e^{\omega (T-T')}w(T)\|^{2}
\]
And thus we have $T_s(T_t x_0) = T_{s+t}x_0$.\\
To see that the semi-group is
\begin{enumerate}
\item a contraction when $\omega > 0$
\item 1-Lipschitz when $\omega = 0$
\item locally Lipschitz in $t$ when $\omega < 0$
\end{enumerate}
we differentiate $\|T_t x_0 - T_t x_1\|^2$ and use equation (23) in conjunction with monotonicity to see that
\[\frac{d}{dt} \|T_t x_0 - T_t x_1\|^2 \leq - \omega \|T_t x_0 - T_t x_1\|^2\]
A simple application of Gronwall's inequality gives the desired conclusions.

\section{Variational resolution of parabolic initial-value problems}

We now apply the results of the last section to the particular class of ASD Lagrangian of the form $L(x,p)=\phi (x)+\phi^*(Ax-p)$ to obtain variational formulations and proofs of existence for  parabolic equations of the form
\begin{eqnarray*}
-\dot x(t) + Ax(t)&\in&\partial \phi (x(t))+\omega x(t)\\
x(0)&=&x_0\\
b_1(x(t)) &=& b_1(x_0).
\end{eqnarray*}
Here again, we have two cases. The first is dealt with in section 5.1 and requires  the operator to be only anti-symmetric while the framework is still purely Hilbertian. The second case  requires that the operator be skew-adjoint --and if necessary-- modulo a pair of boundary operators. The framework there will be on an evolution triple 
$X\subset H \subset X^*$ with $X$ being a Banach space that is anchored on a Hilbert space $H$. It is dealt with in section 5.2.

 \subsection{Parabolic equations involving a diffusion term}
In the  first proposition, we start by assuming the same hypothesis as in Theorem 4.1, that is uniform convexity (in the first variable) of the Lagrangian and a homogeneous initial condition. We will then show how to do away with these conditions in  the corollary that follows.
\begin{proposition}
\label{Asym flow}
 Let $\phi : H\to \bar\R$ be a  convex, lower semi-continuous and proper function on a Hilbert space $H$, and let $A$ be an anti-symmetric linear operator into $H$, with domain $D(A)\subset D(\phi)$.  Assume that:
 \begin{equation}
\hbox{$\phi$ is uniformly convex, with a symmetric domain such that $\partial\phi(0)$ is non-empty.}
 \end{equation}
For any given $\omega\in\R$ and  $T >0$, define the following functional on  $A^2_H([0,T])$
\[
I(u)=\int_0^Te^{2\omega t}\phi(e^{-\omega t}x(t))+e^{2\omega t}\phi^*(e^{-\omega t}(A x(t) - \dot{x}(t)))dt + \frac{1}{2}(\Vert x(0)\Vert^2 + \Vert x(T)\Vert^2).
\]
Then, there exists a path $\bar x \in A^2_H([0,T])$ such that:
\begin{enumerate}
\item
$ I(\bar x)=\inf\limits_{x\in A^2_H([0,T])}I(x)=0$.
\item If $\bar v(t)$ is defined by $\bar v(t) := e^{-\omega t}\bar x(t)$ then it satisfies
\begin{eqnarray}
 \label{eqn:1000}
-\dot{\bar v}(t) + A\bar v(t) -\omega\bar v(t) &\in&\partial\phi(\bar v(t))\quad \mbox{ for a.e. }t\in [0,T] \\
\bar v(0)& = &0. \nonumber
\end{eqnarray}
\end{enumerate}
\end{proposition}
\noindent{\bf Proof:} Setting $\phi_t(x) := e^{2\omega t}\phi(e^{-\omega t}x)$, the  assumptions ensure that
\[L(t, x,p) := \begin{cases}
{\phi_t(x) +\phi_t^*(Ax -p)\mbox{ if }x\in D(\phi)}\cr{+\infty \quad \quad \quad \quad \quad  \mbox{\rm elsewhere}}\cr
\end{cases}\]
is an ASD Lagrangian by Lemma 2.5. Since $\partial\phi(0)$ is  non-empty, it is easy to verify $0\in {\rm Dom}_1\partial L$ and all the hypothesis of Theorem 4.1 are satisfied. Therefore, there exists $\bar x(\cdot)\in A^2_H([0,T])$ such that
\[ 0 = \int_0^T L(t,\bar x(t),\dot {\bar x}(t))\, dt+\ell\big( \bar x(0),\bar x(T)\big)\]
 $\ell(a,b) := \frac{1}{2}(\Vert a\Vert^2 + \Vert b\Vert^2)$. 
Therefore
\begin{eqnarray*}
0&=& \int_0^T\phi_t(\bar x(t))+\phi_t^*(A\bar x(t) - \dot{\bar x}(t))dt + \frac{1}{2}(\Vert\bar x(0)\Vert^2 + \Vert\bar x(T)\Vert^2)\\
&\geq& \int_0^T\langle \bar x(t),A\bar x(t) - \dot{\bar x}(t)\rangle dt + \frac{1}{2}(\Vert\bar x(0)\Vert^2 + \Vert\bar x(T)\Vert^2)\\
&=&\Vert\bar x(0)\Vert^2\geq 0.
\end{eqnarray*}
It follows that,
\[-\dot{\bar x}(t) + A\bar x(t)\in e^{\omega t}\partial\phi(e^{-\omega t}\bar x(t))\]
\[\bar x(0) =0\]
and by a simple application of the product-rule  we see that $\bar v(t)$ defined by $\bar v(t) := e^{-\omega t}\bar x(t)$ satisfies (\ref{eqn:1000}).
 
\begin{corollary}
\label{asym flow with initial condition}
Let $\phi : H\to \bar\R$ be a  convex, lower semi-continuous and proper function on a Hilbert space $H$, and let $A$ be an anti-symmetric linear operator into $H$, with domain $D(A)\subset D(\phi)$.  Assume that:
\begin{equation}
 x_0\in D(A)\cap \partial\phi.
 \end{equation}
 Then, for all $\omega\in\R$ and for all $T >0$, there exists $\bar u\in A^2_H([0,T])$ such that
\begin{eqnarray}
\label{Eq:1001}
-\dot{\bar u}(t) + A\bar u(t) -\omega\bar u(t) &\in&\partial\phi(\bar u(t))\quad \mbox{ for a.e. }t\in [0,T]\\
\bar u(0)& = &x_0.\nonumber
\end{eqnarray}
\end{corollary}
\noindent{\bf Proof:}
Define the convex function $\psi : H\to\bar\R$ by
\[
\psi(x) := \phi(x + x_0) + \frac{\Vert x\Vert^2}{2} -\langle x,Ax_0\rangle + \langle x,\omega x_0\rangle.
\]
By the fact that $\partial\phi(x_0)$ is non-empty, it is easy to check that $\psi$ satisfies all the conditions of Proposition $\ref{Asym flow}$. Therefore, there exists $\bar v(\cdot)\in A^2_H([0,T])$ satisfying the evolution equation
\[-\dot{\bar v}(t) + A\bar v(t) \in\partial\psi(\bar v(t))+(\omega-1)\bar v(t)\mbox{ for a.e. }t\in [0,T]\]
\[\bar v(0) = 0\]
Since $\partial\psi(x) = \partial\phi(x + x_0) + x -Ax_0 + \omega x_0$, we get that $\bar u(t) := \bar v(t) + x_0$ satisfies equation (\ref{Eq:1001}).

\subsubsection*{Evolution driven by the transport operator and the $p$-Laplacian}
Consider the following evolution equation on a smooth bounded domain of $\R^n$. 
\begin{eqnarray}
\label{Eq:transport.plus.heat}
- u_t(x,t)+ \vec a(x)\cdot\nabla  u(x,t) &=& -\Delta_p u(x,t) + \frac{1}{2}a_0(x) u(x,t)+\omega u(x,t) \quad \hbox{\rm on $[0,T]\times \Omega$}\\
 u(x,0)& =& u_0(x)\quad \hbox{\rm on $\Omega$}\nonumber \\
 u(x,t) &=& 0\quad \quad \quad \hbox{\rm on $[0,T]\times \partial \Omega$.} \nonumber
\end{eqnarray}
We can establish variationally the following
\begin{corollary}
Let $\vec a:\R^n\to\R^n$ be a smooth vector field and $a_0\in L^\infty(\Omega)$. 
For $p\geq 2$, $\omega\in\R$, and any $u_0$  in $W^{1,p}_0(\Omega)\cap \lbrace u; \Delta_p u\in L^2(\Omega)\rbrace$,  there exists $\bar u \in A^2_{L^2(\Omega)}([0,T])$ that solves (\ref{Eq:transport.plus.heat}).
Furthermore, $\Delta_p\bar u(x,t)\in L^2(\Omega)$ for almost all $t\in[0,T]$.
\end{corollary}
\noindent{\bf Proof:} The operator $Au = \vec a\cdot\nabla u + \frac{1}{2}(\nabla\cdot\vec a)u$ 
with domain $D(A) = H^1_0(\Omega)$ is anti-symmetric.  In order to apply corollary $\ref{asym flow with initial condition}$ with $H = L^2(\Omega)$ and $A$,  we need to insure convexity of the potential and for that we pick $K >0$ such that $\nabla\cdot\vec a(x) + a_0(x) + K \geq 1$ for all $x\in\Omega$.\\
 Now define $\phi : H \to \bar\R$ by
\[\phi(u) := 
\begin{cases}{\frac{1}{p}\int_\Omega \vert\nabla u(x)\vert^p dx + \frac{1}{4}\int_\Omega(\nabla\cdot\vec a(x) + a_0(x) + K)\vert u(x)\vert^2dx\quad \quad {\rm if} \quad  u\in W^{1,p}_0(\Omega)}\cr{+\infty \quad \quad \quad \quad {\rm elsewhere}}
\end{cases}\]
By observing that $\phi$ is a convex l.s.c. function with symmetric domain and $D(\phi)\subset D(A)$, we can apply Corollary $\ref{asym flow with initial condition}$ with the linear factor ($\omega-\frac{K}{2}$),   to obtain the existence of a $\bar u(\cdot)\in A^2_H([0,T])$ such that
\[-\dot{\bar u}(t) + A\bar u(t) \in\partial\phi(\bar u(t))+(\omega-\frac{K}{2})\bar u(t)\mbox{ for a.e. }t\in [0,T]\]
\[\bar u(0) = x_0\]
and this is precisely the equation (\ref{Eq:transport.plus.heat}).
Since now $\partial\phi(\bar u(t))$ is a non-empty set in $H$ for almost all $ t\in[0,T]$, we have $\Delta_p \bar u(x,t)\in L^2(\Omega)$ for almost all $t\in[0,T]$.  

\subsection{Parabolic equations driven by first-order operators}
In  this subsection we deal with parabolic equations of the form:
\begin{eqnarray}
\label{Eqn:1002}
-\dot x(t) + Ax(t)+w x(t)&\in&\partial \phi (x(t))+w x(t)\quad \mbox{\rm for a.e. }t\in [0,T]\nonumber\\
x(0)&=&x_0\\
b_1(x(t)) &=& b_1(x_0). \quad \mbox{\rm for a.e. }t\in [0,T].\nonumber
\end{eqnarray}
where the operator $A$ is skew-adjoint modulo  boundary operators $(b_1,b_2)$. Here we need the framework of an evolution triple, where  $X$ is a reflexive Banach space and $H$ is a Hilbert space satisfying $X \subset H \subset X^*$. in such a way that each space is dense in the following one. Again we start with a theorem that assumes all the hypothesis of Theorem 4.1. We will then relax these conditions in the corollary that follows it. 

\begin{proposition} \label{ASAMB flow}
Let $X\subset H \subset X^*$ be an evolution triple and let $A:D(A)\subset X\to X^*$ be a skew-adjoint operator modulo boundary operators $(b_1, b_2): D(b_1,b_2) \to H$. Let $\phi :X\to\R$ be a  convex lower semi-continuous and proper function on $X$,  that is bounded on the bounded sets of $X$ and also coercive on $X$. Assume that 
\begin{equation}
\hbox{$\phi$ is uniformly convex on $H$ and  $\partial\phi(0)\cap H$ is non-empty. }
\end{equation}
 Let $w\in\R$ and $T>0$, then there exists a solution $v \in A_H^2$ for the initial value problem
\begin{eqnarray}
 \label{initialvalue}
 -\dot v(t)+Av(t)&\in& \partial \phi \big( v(t)\big) + wv(t)\quad \mbox{\rm a.e.} \quad t\in [0,T]\\
 b_1(v(t)) &=& 0 \quad \mbox{\rm for a.e. }t\in [0,T]  \nonumber\\
 v(0)&=&0. \nonumber
\end{eqnarray}
 It is obtained by minimizing over  $A_H^2$ the functional   
 \begin{eqnarray*}
 I(u)&=&\int_0^Te^{2\omega t}\left\{\phi (e^{-\omega t} u(t))+\phi^*(e^{-\omega t}(-A^au(t) -{\dot u(t))}\right\}dt\\
 && + \frac{1}{2}\int_0^T(\Vert b_1(x(t))\Vert^2_{H_1} + \Vert b_2(x(t))\Vert^2_{H_2})dt  +\frac{1}{2}\|u(0)\|^{2} +\frac{1}{2}\|u(T)\|^{2}.
 \end{eqnarray*}
The minimum of $I$ is then zero and is attained at a path $y(t)$ such that  $x(t)= e^{-\omega t} y(t)$ is a solution of (\ref{initialvalue}). 
\end{proposition}
 Typical convex functions satisfying the conditions above are ones such that for some $C>0$, $m,n>1$ we have the following growth condition:
\begin{equation}
\label{growth}
\C\left( {\| x\|}_X^m-1\right)\le\phi (x)\le C\left( {\| x\|}_{X}^n+1\right).
\end{equation}
The corresponding Lagrangian $L$ is ASD on $X\times X^*$ where $X\subseteq H\subseteq X^*$.  Since our theory for evolution equations applies to Hilbert spaces, the following lemma will bridge the gap:
\begin{lemma}\label{lift}
Let $X\subset H \subset X^*$ be an evolution triple, and suppose $L:X\times X^*\to\R$ is ASD on the Banach space $X$.  Assume the following two conditions:
\begin{enumerate}
\item For all $x\in X$, the map $L(x,\cdot ):X^*\to\overline{\R}$ is  continuous on $X^*$.
\item There exists $x_0\in X$ such that $p\to L(x_0,p)$ is bounded on the bounded sets of $X^*$.
\end{enumerate} 
Then the Lagrangian defined on $H$ by 
\begin{eqnarray*}
M(x,p):=\left\{ \begin{array}{ll}L(x,p) &x\in X\\
+\infty &x\in H\backslash X\end{array}\right.
\end{eqnarray*}
is anti-selfdual on $H\times H$.
\end{lemma}
 {\bf Proof:} For $(\tilde x,\tilde p)\in X\times H$, write
\begin{eqnarray*}
M^*(-\tilde p,-\tilde x) &=& \sup_{\stackrel{x\in X}{p\in H}}\left\{ {\braket{\tilde x}{p}}_H 
  +{\braket{\tilde p}{x}}_H-L(x,p)\right\}\\
&=&\sup_{X\in X}\sup_{p\in H}\left\{ \phantom{x}_X{\braket{\tilde x}{p}}_{X^*}
  +\phantom{x}_X{\braket{x}{\tilde p}}_{X^*}
   -L(x,p)\right\}\\
&=&\sup_{X\in X}\sup_{p\in X^*}\left\{ \phantom{x}_x{\braket{\tilde x}{p}}_{X^*}
  +\phantom{x}_X{\braket{x}{\tilde p}}_{X^*}
   -L(x,p)\right\}\\
&=& L(-\tilde x,-\tilde p)
\end{eqnarray*}
Now suppose $\tilde x\in H\backslash X$. Then
\begin{eqnarray*}
M^*(\tilde p,\tilde x) &=& \sup_{\stackrel{x\in X}{p\in H}} \left\{ {\braket{\tilde x}{p}}_H
   +{\braket{\tilde p}{x}}_H
   -L(x,p)\right\}\\
&\ge& \braket{\tilde p}{x_0} +\sup_{p\in H}\left\{ {\braket{\tilde x }{p}}_H
   -L(x_0,p)\right\} 
\end{eqnarray*}
Since $\tilde x\notin X$,  we have that $\sup \left\{ \braket{\tilde x }{p}; p\in H, \|p\|_{X^*}\leq 1\right\} =+\infty $. Since $p\to L(x_0,p)$ is bounded on the bounded sets of $X^*$, it follows that 
\[
M^*(\tilde p,\tilde x) \geq  \braket{\tilde p}{x_0} +\sup_{p\in H}\left\{ {\braket{\tilde x }{p}}_H
   -L(x_0,p)\right\}=+\infty,
\] 
and we are done. \\

\noindent{\bf Proof of Proposition 5.2:} Again
\[L(x,p) :=
\begin{cases}
{\phi(x) +\phi^*(Ax - p) + \frac{1}{2}(\Vert b_1(x)\Vert^2_{H_1} + \Vert b_2(x)\Vert^2_{H_2})\quad {\rm if} \ x\in D(A)\cap D(b_1, b_2)}\cr {+\infty \quad \quad \quad \quad \quad \quad \quad {\rm elsewhere}}\cr
\end{cases} \]
is an ASD Lagrangian on $X\times X$ by Propostion 2.1. The coercivity condition on $\phi$ ensures --via Lemma \ref{lift}--that $L(x,p)$ lifts to a ASD Lagrangian on $H\times H$ that is uniformly convex in the first variable. It is easy to check that all the conditions of Theorem 4.1 are satisfied by $L(x,p)$. Therefore, there exists $\bar x(\cdot)\in A^2_H([0,T])$ such that $I(\bar u)=0$, which yields
\begin{eqnarray*}
0&=& \int_0^T\phi_t(\bar x(t))+\phi_t^*(A\bar x(t) - \dot{\bar x}(t)) + \frac{1}{2}(\Vert b_1(x(t))\Vert^2_{H_1} + \Vert b_2(x(t))\Vert^2_{H_2})dt + \frac{1}{2}(\Vert\bar x(0)\Vert^2 + \Vert\bar x(T)\Vert^2)\\
&\geq& \int_0^T\langle \bar x(t),A\bar x(t) - \dot{\bar x}(t)+ \frac{1}{2}(\Vert b_1(x(t))\Vert^2_{H_1} + \Vert b_2(x(t))\Vert^2_{H_2})dt + \frac{1}{2}(\Vert\bar x(0)\Vert^2 + \Vert\bar x(T)\Vert^2)\\
&=&\int_0^T\Vert b_1(\bar x(t))\Vert^2_{H_1}dt + \Vert\bar x(0)\Vert^2\geq0.
\end{eqnarray*}
So all inequalities are equalities, and we obtain $-\dot{\bar x}(t) + A\bar x(t)\in e^{\omega t}\partial\phi(e^{-\omega t}\bar x(t))$, $b_1(\bar x(t)) = 0$, and $\bar x(0) =0$. We now set   $\bar v(t) := e^{-\omega t}\bar x(t)$ and the rest is straightforward.
 \begin{corollary}
\label{asamb flow with initial condition}
Let $X\subset H \subset X^*$ be an evolution triple and let $A:D(A)\subset H\to X^*$ be a skew-adjoint operator modulo boundary operators $(b_1, b_2): D(b_1,b_2) \to H$. Let $\phi :X\to\R$ be a  convex lower semi-continuous and proper function on $X$,  that is bounded on the bounded sets of $X$ and also coercive on $X$. Assume that 
\begin{equation}
\hbox{$x_0\in D(A)\cap D(b_1, b_2)$ and $ \partial\phi(x_0)\cap H$ is non-empty.}
\end{equation}
 Then, for all $\omega\in\R$ and for all $T >0$, there exists $\bar u \in A^2_H([0,T])$ which solves (\ref{Eqn:1002}).
 \end{corollary}
\noindent{\bf Proof:}
Define the convex function $\psi : X\to\R$ by
\[\psi(x) := \phi(x + x_0) + \frac{\Vert x\Vert^2}{2} -\langle x,Ax_0\rangle + \langle x,\omega x_0\rangle\]
It is easy to check that $\psi$ satisfies all the conditions of Proposition $\ref{ASAMB flow}$. Therefore, there exists $\bar v \in A^2_H([0,T])$ satisfying the evolution equation
\[-\dot{\bar v}(t) + A\bar v(t) \in\partial\psi(\bar v(t))+(\omega-1)\bar v(t)\mbox{ for a.e. }t\in [0,T]\]
\[b_1(\bar v(t)) = 0\mbox{ for a.e. }t\in [0,T]\]
\[\bar v(0) = 0.\]
Since $\partial\psi(x) = \partial\phi(x + x_0) + x -Ax_0 + \omega x_0$, we have that $\bar u(t) := \bar v(t) + x_0$  solves (\ref{Eqn:1002}).

\subsubsection*{Evolutions driven by transport operators}

Consider the evolution equation
\begin{eqnarray}
\label{Eq:transport.No.heat}
-u_t(x,t) + \vec a(x)\cdot\nabla u(x,t) &=& \frac{1}{2} a_0(x)u(x,t) + u(x,t)\vert u(x,t)\vert^{p-2} + \omega u(x,t)\quad \hbox{\rm on $[0,T]\times \Omega$}\nonumber\\
u(x,0) &=& u_0(x)\quad \hbox{\rm on $\Omega$} \\
u(x,t) &= &u_0(x) \quad \hbox{\rm on $[0,T]\times \Sigma_+$}\nonumber
\end{eqnarray}
We assume that the domain $\Omega$ and the vector field $\vec a(\cdot)$ satisfies all the assumption in section 1.3.

\begin{corollary}
Let $p > 1$, $f\in L^2(\Omega)$ and $a_0\in L^{\infty}(\Omega)$. For any $\omega \in\R$ and $u_0 \in L^{\infty}(\Omega)\cap H^1(\Omega)$ there exists $\bar u(\cdot)\in A^2_{L^2(\Omega)}([0,T])$ satisfying (\ref{Eq:transport.No.heat}).
\end{corollary}
\noindent{\bf Proof}: We distinguish two cases:\\

\noindent {\bf Case 1: $p\geq 2$}. We then take $X = L^p(\Omega)$, $H = L^2(\Omega)$. since again the operator $A : D(A)\to X^*$ defined as $Au=\vec{a}\cdot\nabla u+\frac{1}{2}\left( \nabla\cdot\vec{a}\right) u$  with domain
\[
D(A) = \lbrace u\in L^p(\Omega) \mid \vec a\cdot\nabla u+\frac{\nabla\vec a}{2} u\in L^q(\Omega)\rbrace
\]
is skew-adjoint modulo the boundary operators $(b_1u, b_2u) = (u\vert_{\Sigma_+},u\vert_{\Sigma_-})$
whose domain  is 
\[
D(b_1,b_2) = \lbrace u\in L^p(\Omega) \vert (u\vert_{\Sigma_+},u\vert_{\Sigma_-})\in L^2(\Sigma_+ ; \vert \vec a \cdot\hat n\vert d\sigma)\times L^2(\Sigma_- ; \vert \vec a \cdot\hat n \vert d\sigma)\rbrace
\]
 
\noindent {\bf Case 2: $1<p<2$}.  The space is then $X = H = L^2(\Omega)$. \\

In both case, we pick $K >0$ such that $\nabla\cdot\vec a(x) + a_0(x) + K \geq 1$ for all $x\in\Omega$, and define the function $\phi : X\to \R$ by
\[\phi(u):=\frac{1}{p}\int_\Omega \vert u(x)\vert^p dx+ \frac{1}{4}\int_\Omega(\nabla\cdot\vec a(x) + a_0(x) + K)\vert u(x)\vert^2dx\]
Then $\phi$ is a convex, l.s.c. function that is bounded on bounded sets of $X$ and coercive on $X$. Since $u_0(x)\in L^{\infty}(\Omega)\cap H^1(\Omega)$, $\partial\phi(u_0)$ is non-empty and $u_0 \in D(A)\cap D(b_1, b_2)$. 

So by corollary $\ref{asamb flow with initial condition}$, there exists $\bar u(\cdot)\in A^2_H([0,T])$ such that
\[-\dot{\bar u}(t)+A \bar u(t)\in\partial\phi(\bar u(t))+(\omega -\frac{K}{2})\bar u(t)\ quad {\rm for}\quad t\in [0,T]\]
\[b_1(\bar u(t))= b_1(u_0)\quad {\rm for}\quad  t\in [0,T]\]
\[\bar u(0) = u_0\]
and this is precisely the equation (\ref{Eq:transport.No.heat})
and we are done.

\end{document}